\documentclass[aos,seceqn,MSNbibl,citesort,dvips]{arximspdf}
\usepackage{graphicx}

\doi{10.1214/10-AOS864}
\volume{39}
\issue{2}
\pubyear{2011}
\firstpage{1180}
\lastpage{1210}

\makeatletter

\def\@bmisc[#1]{%
  \get@battribute{unstr}%
  \common@pub@types%
  \let\bauthor\bbl@bauthor%
  \let\bhowpublished\@firstofone%
  \def\borganization##1{{\bauthor@style ##1}}%
}

\def\bptnote#1{}

\newtheorem{lemma}{Lemma}[section]
\newtheorem{theorem}{Theorem}[section]
\newtheorem{corollary}{Corollary}[section]
\newproclaim{assumption}{Assumption}[section]
\newproclaim{exam}{A toy example}

\newproclaim{Remarks}{Remarks}
\newproclaim{Remark}{Remark}

\makeatother

\begin{document}
\begin{frontmatter}

\title{Performance guarantees for individualized treatment rules\thanksref{T1}}
\runtitle{Individualized treatment rules}

\thankstext{T1}{Supported by NIH Grants R01 MH080015 and P50 DA10075.}

\begin{aug}
\author[A]{\fnms{Min} \snm{Qian}\corref{}\ead[label=e1]{minqian@umich.edu}}
and
\author[A]{\fnms{Susan A.} \snm{Murphy}\ead[label=e2]{samurphy@umich.edu}}
\runauthor{M. Qian and S. A. Murphy}
\affiliation{University of Michigan}
\address[A]{Department of Statistics\\
University of Michigan\\
439 West Hall\\
1085 S. University Ave.\\
Ann Arbor, Michigan 48109\\
USA\\
\printead{e1}\\
\phantom{E-mail: }\printead*{e2}} %adresu isvedimo komanda gale!
\end{aug}

% HISTORY:
\received{\smonth{1} \syear{2010}}
\revised{\smonth{11} \syear{2010}}

% ABSTRACT
%
\begin{abstract}
Because many illnesses show heterogeneous response to treatment, there
is increasing interest in individualizing treatment to patients
[\textit{Arch. Gen. Psychiatry} \textbf{66} (2009) 128--133].
An \textit{individualized treatment rule} is a
decision rule that recommends treatment according to patient
characteristics. We consider the use of clinical trial data in the
construction of an individualized treatment rule leading to highest
mean response. This is a difficult computational problem because the
objective function is the expectation of a weighted indicator function
that is nonconcave in the parameters. Furthermore, there are
frequently many pretreatment variables that may or may not be useful in
constructing an optimal individualized treatment rule, yet cost and
interpretability considerations imply that only a few variables should
be used by the individualized treatment rule. To address these
challenges, we consider estimation based on $l_1$-penalized least
squares. This approach is justified via a finite sample upper bound on
the difference between the mean response due to the estimated
individualized treatment rule and the mean response due to the optimal
individualized treatment rule.
\end{abstract}

% KEYWORDS
%

\end{frontmatter}

%s1 ###
\section{Introduction}\label{sec1}

Many illnesses show heterogeneous response to treatment. For
example, a study on schizophrenia \cite{ishigooka2001} found that
patients who take the same antipsychotic (olanzapine) may have very
different responses. Some may have to discontinue the treatment due
to serious adverse events and/or acutely worsened symptoms, while
others may experience few if any adverse events and have improved
clinical outcomes. Results of this type have motivated
researchers to advocate the individualization of treatment to each
patient \cite{lesko2007,piquette2007,insel2009}. One step in this
direction is to estimate each patient's risk level and then match
treatment to risk category \cite{cai2008,cai2010}. However, this
approach is best used to decide whether to treat; otherwise it
assumes the knowledge of the best treatment for each risk category.
%Early research in this direction focused on identifying possible
%qualitative interactions between treatment and patient pretreatment
%variables via hypothesis testing and then the best treatment for
%different subset of patients is selected according to the test
%results
% (\citealt{byar1977,shuster1983,gail1985}).
Alternately, there is an abundance of literature
focusing on predicting each patient's prognosis
under a particular treatment \cite{feldstein1978,stoehlmacher2004}.
Thus, an obvious way to individualize treatment is to recommend the
treatment achieving the best predicted prognosis for that patient. In
general, the goal is to use data to construct individualized treatment rules
that, if implemented in future, will optimize
the mean response.

%When data in which patients are randomized to one of several active
%treatments are available, an improved variant is to use a prediction
%model that involves the treatment indicator as well as interactions
%between the treatment indicator and covariates.
%In this paper,
%we consider an improved variant of this ``prognosis prediction
%approach." We directly estimate an individualized treatment rule that
%recommends treatment according to individual characteristics and we
%provide a finite sample performance bound on the quality of the
%resulting rule.

Consider data from a single stage randomized trial involving
several active treatments. A first natural procedure to construct the
optimal individualized treatment rule is to maximize an empirical
version of the mean response over a class of treatment rules (assuming
larger responses are preferred). As will be seen, this maximization
is computationally difficult because the mean response of a
treatment rule is the expectation of a weighted indicator that is
noncontinuous and nonconcave in the parameters. To address this
challenge, we make a substitution. That is, instead of directly maximizing
the empirical mean response to estimate the treatment rule, we
use a two-step procedure that first estimates a
conditional mean and then from this estimated conditional mean derives the
estimated treatment rule. As will be seen in Section \ref
{sec:relation}, even if the optimal treatment rule is contained in the
space of treatment rules considered by the substitute two-step
procedure, the estimator derived from the two-step procedure may not be
consistent.
%; this is due to the mismatch between maximizing the mean response
%over treatment rules and the two-stage procedure.
However, if the conditional mean is modeled correctly, then
the two-step procedure consistently estimates the optimal
individualized treatment rule.
This
motivates consideration of rich conditional mean models with many
unknown parameters.
Furthermore, there
are frequently many pretreatment variables that may or may not be
useful in constructing an optimal individualized treatment rule, yet
cost and interpretability considerations imply that
fewer rather than more variables should be used by the treatment
rule. This consideration motivates the use of $l_1$-penalized least
squares ($l_1$-PLS).
%, it is problematic to assume that a
%correct or nearly correct model is used to analyze the data.
%Of course our
%finite sample performance bounds are tighter the more closely our
%model approximates the truth.
%The
%tuning parameter in the $l_1$ penalty can be selected by maximizing
%the estimated mean response. As compared to treatment assignment via
%separate prognosis prediction for each treatment, our treatment rule
%will employ much fewer patient pretreatment variables and will yield
%higher mean response when the approximation model is poor.

We propose to estimate an optimal individualized treatment rule using
a~two step procedure that first estimates the conditional mean response
using $l_1$-PLS with a rich linear model and second, derives the
estimated treatment rule from estimated conditional mean.
For brevity, throughout, we call the two step procedure the $l_1$-PLS method.
We
derive several finite sample upper bounds on the difference between
the mean response to the optimal treatment rule and the mean
response to the estimated treatment rule. All of the
upper bounds hold even if our linear model for the conditional mean
response is incorrect and
to our knowledge are, up to constants, the best available.
We use the upper bounds in Section~\ref{sec:relation} to illuminate the
potential
mismatch between using least squares in the two-step procedure and the
goal of
maximizing the mean response.
The
upper bounds in Section~\ref{sec:finaloracle} involve a minimized sum
of the approximation error and estimation
error; both errors result from the estimation of the conditional mean response.
%The estimation error is smaller if the linear model is sparse with few
%parameters or if the sample size is large.
We shall see that $l_1$-PLS estimates a linear model that minimizes
this approximation plus estimation error sum among a set of suitably
sparse linear models.

If the part of the model for the conditional mean
involving the treatment effect is correct, then the upper bounds imply
that, although a surrogate two-step procedure is used, the estimated
treatment rule is consistent. The upper bounds provide a convergence
rate as well. Furthermore, in this
setting, the upper bounds can be used to inform how
to choose the tuning parameter involved in the $l_1$ penalty to
achieve the best rate of convergence. As a~by-product,
this paper also contributes to existing literature on $l_1$-PLS
by providing a finite sample prediction error
bound for the $l_1$-PLS estimator in the random design setting without
assuming the model class contains or is close to the true model.

%Lastly we use the upper bound to intuitively illustrate the connection
%between $l_1$-penalization and
%$l_0$-penalization in prediction.

The paper is organized as follows. In Section \ref{sec:prelim}, we
formulate the decision making problem. In Section
\ref{sec:relation}, for any given decision, that is, individualized
treatment rule, we relate the reduction in mean response to the excess
prediction error. In Section \ref{sec:lasso}, we estimate an optimal
individualized treatment rule via
$l_1$-PLS and provide a finite sample upper
bound on the reduction in mean response achieved by the estimated rule.
In Section \ref{sec:data}, we consider a data dependent tuning
parameter selection criterion. This method is evaluated using
simulation studies and illustrated with data from the
Nefazodone-CBASP trial \cite{keller2000}. Discussions and future
work are presented in Section \ref{sec:discussion}.

%s2 ###
\section{Individualized treatment rules} \label{sec:prelim}
We use upper case letters to denote random variables and lower case
letters to denote values of the random variables. Consider data from
a randomized trial. On each subject, we have the pretreatment
variables $X\in\mathcal{X}$, treatment $A$ taking values in a
finite, discrete treatment space~$\mathcal{A}$, and a real-valued
response $R$ (assuming large values are desirable). An
\textit{individualized treatment rule} (ITR) $d$ is a deterministic
decision rule from $\mathcal{X}$ into the treatment space
$\mathcal{A}$.

Denote the distribution of $(X,A,R)$ by $P$. This is the
distribution of the clinical trial data; in particular, denote the
known randomization distribution of $A$ given $X$ by $p(\cdot|X)$.
The likelihood of
$(X,A,R)$ under $P$ is then $f_0(x)p(a|x)f_1(r|x,a)$, where $f_0$ is
the unknown density of
$X$ and $f_1$ is the unknown density of $R$ conditional on $(X,A)$.
Denote the expectations with respect to the distribution $P$ by an
$E$.
For any ITR $d\dvtx
\mathcal{X}\rightarrow\mathcal{A}$, let $P^d$ denote the
distribution of $(X,A,R)$ in which $d$ is used to assign treatments.
Then the likelihood of $(X,A,R)$ under $P^{d}$ is
$f_0(x)1_{a=d(x)}f_1(r|x,a)$. Denote
expectations with respect to the distribution $P^d$ by an $E^d$. The
\textit{Value} of~$d$ is defined as $V(d)\triangleq E^d(R)$. An \textit{optimal
ITR}, $d_0$, is a rule that has the maximal Value,
that is,
\[
d_0\in\mathop{\arg\max}_d V(d),
\]
where the $\arg\max$ is over all possible decision rules. The Value of
$d_0$, $V(d_0)$, is the \textit{optimal Value}. %Note there may
%be multiple decision rules that maximize the value.

Assume $P[p(a|X)>0]=1$ for all $a\in\mathcal{A}$ (i.e., all
treatments in $\mathcal{A}$ are possible for all values of $X$
a.s.). Then $P^d$ is absolutely continuous with respect to $P$ and a
version of the Radon--Nikodym derivative is
$dP^d/dP=1_{a=d(x)}/p(a|x)$. Thus, the Value of $d$ satisfies
%
%e2.1 ###
\begin{equation}\label{eqn:value}
V(d)=E^d(R)=\int R\,dP^d=\int
R\,\frac{dP^d}{dP}\,dP=E\biggl[\frac{1_{A=d(X)}}{p(A|X)}R\biggr].
\end{equation}
Our goal is to estimate $d_0$, that is, the ITR that
maximizes (\ref{eqn:value}), using data from distribution $P$. When
$X$ is low dimensional and the best rule within a~simple class of
ITRs is desired, empirical versions of the Value can be
used to construct estimators \cite{murphy2001,robins2008}. However,
if the best rule within a larger class of ITRs is of
interest, these approaches are no longer feasible.

Define $Q_0(X,A)\triangleq E(R|X,A)$ [$Q_0(x,a)$ is sometimes called
the ``Quality'' of treatment $a$ at observation $x$]. It follows from
(\ref{eqn:value}) that
for any ITR $d$,
\[
V(d) =
E\biggl[\frac{1_{A=d(X)}}{p(A|X)}Q_0(X,A)\biggr]=E\biggl[\sum_{a\in\mathcal
{A}}1_{d(X)=a}Q_0(X,a)\biggr]
=E[Q_0(X,d(X))].
\]
Thus, $V(d_0)=E[Q_0(X,d_0(X))]\leq E[\max_{a\in\mathcal{A}} Q_0(X,a)]$.
On the other hand,
by the definition of $d_0$,
\[
V(d_0)\geq
V(d)|_{d(X)\in\mathop{\arg\max}_{a\in\mathcal{A}}
Q_0(X,a)}=E\Bigl[\max_{a\in\mathcal{A}} Q_0(X,a)\Bigr].
\]
Hence, an optimal ITR
satisfies $d_0(X)\in\arg\max_{a\in\mathcal{A}}$ $Q_0(X,a)$ a.s.

%s3 ###
\section{Relating the reduction in Value to excess prediction error}
\label{sec:relation}

The above argument indicates that the estimated ITR will be of high quality
(i.e., have high Value) if we can estimate $Q_0$
accurately. In this section, we justify this by providing a
quantitative relationship between the Value and
the prediction error.

Because $\mathcal A$ is a finite, discrete treatment space, given any ITR,
$d$, there \mbox{exists} a square integrable function
$Q\dvtx\mathcal{X}\times\mathcal{A}\rightarrow\mathbb{R}$ for which
$d(X)\in\break\arg\max_aQ(X,a)$ a.s. Let $L(Q)\triangleq E[R-Q(X,A)]^2$
denote the prediction error of $Q$ (also called the mean quadratic loss).
Suppose that $Q_0$ is square integrable and that the randomization
probability satisfies $p(a|x)\geq
S^{-1}$ for an $S>0$ and all $(x,a)$ pairs. Murphy \cite{murphy2005}
showed that
%
%e3.1 ###
\begin{equation} \label{eqn:bound1}
V(d_0)-V(d)\leq
2S^{1/2}[L(Q)-L(Q_0)]^{1/2}.
\end{equation}
Intuitively, this upper bound means that if the excess prediction error of
$Q$ [i.e., $L(Q)-L(Q_0)$] is small, then the reduction in Value of the
associated ITR $d$ [i.e., $V(d_0)-V(d)$] is small.
Furthermore, the upper bound provides a
rate of convergence for the Value of an estimated ITR. For example, suppose
$Q_0$ is linear, that is, $Q_0=\Phi(X,A)\bolds{\theta}_0$ for a
given vector-valued basis function $\Phi$ on $\mathcal{X}\times\mathcal
{A}$ and an unknown parameter ${\theta}_0$.
And suppose we use a correct linear model for $Q_0$ (here ``linear''
means linear in parameters), say the model $\mathcal{Q}=\{\Phi
(X,A)\bolds{\theta}\dvtx\bolds{\theta}\in\mathbb{R}^{\mathrm{dim}(\Phi)}\}
$ or a linear model containing $\mathcal{Q}$ with dimension of
parameters fixed in $n$. If we estimate $\bolds{\theta}$ by least
squares and denote the estimator by $\bolds{\hat\theta}$, then the
prediction error of $\hat Q =\Phi\bolds{\hat\theta}$ converges to
$L(Q_0)$ at rate $1/n$ under mild regularity conditions.
This together with inequality (\ref{eqn:bound1}) implies that
the Value obtained by the estimated ITR, $\hat d(X)\in\arg\max_a\hat
Q(X,a)$, will converge to
the optimal Value at rate at least $1/\sqrt n$.

In the following theorem, we improve this upper
bound in two aspects. First, we show that an upper bound with
exponent larger than $1/2$ can be obtained under a margin condition,
which implicitly implies a faster rate of convergence.
Second, it turns out that the upper bound need only depend on one term
in the function $Q$; we call this the treatment effect term $T$. For
any square integrable $Q$, the associated treatment effect
term is defined as $T(X,A)\triangleq Q(X,A) - E[Q(X,A)|X]$. Note that
$d(X)\in\arg\max_a T(X,a)=\arg\max_a Q(X,a)$ a.s.
Similarly, the true treatment effect
term is given by
%
%e3.2 ###
\begin{equation} \label{eqn:trteffect}
T_0(X,A)\triangleq Q_0(X,A) - E[Q_0(X,A)|X].
\end{equation}
$T_0(x,a)$ is the centered effect of treatment
$A=a$ at observation $X=x$; $d_0(X)\in\arg\max_a T_0(X,a)$.
%The following result implies that if $T$ is close to $T_0$ then the
%value of $d$ is close to the value of $d_0$.
%Under appropriate conditions, consistency of estimators of
%$T_0$ will be unaffected by the whether the estimator of
%$E[Q_0(X,A)|X]$ is consistent. See
%Section \ref{sec:lasso} for further discussion.
%
%In classification tighter upper bounds can be obtained
%under margin conditions (\citealt{bartlett2006}).
%
\begin{theorem}\label{thm:bound2}
Suppose $p(a|x)\geq S^{-1}$ for a positive constant $S$ for all
$(x,a)$ pairs. Assume there exists some constants $C>0$ and
$\alpha\geq0$ such that
%
%e3.3 ###
\begin{equation}\label{eqn:noise}
\mathbf{P}\Bigl(\max_{a\in\mathcal{A}}T_0(X,a)-
\max_{a\in\mathcal{A}\setminus\mathop{\arg\max}_{a\in\mathcal{A}}T_0(X,a)}T_0(X,a)
\leq\epsilon\Bigr)\leq C\epsilon^\alpha
\end{equation}
for all positive $\epsilon$.
Then for any ITR
$d\dvtx\mathcal{X}\rightarrow\mathcal{A}$ and square integrable function
$Q\dvtx\mathcal{X}\times\mathcal{A}\rightarrow\mathbb{R}$ such that
$d(X)\in\arg\max_{a\in\mathcal{A}}Q(X,a)$ a.s., we have
%
%e3.4 ###
\begin{equation}
\label{eqn:bound2}
V(d_0)-V(d)\leq C' [L(Q)-L(Q_0)]^{(1+\alpha)/(2+\alpha)}
\end{equation}
and
%
%e3.5 ###
\begin{equation}
\label{eqn:bound3}
V(d_0)-V(d)\leq C'
\bigl[E\bigl(T(X,A)-T_0(X,A)\bigr)^2\bigr]^{(1+\alpha)/(2+\alpha)},
\end{equation}
where $C'=(2^{2+3\alpha}S^{1+\alpha}C)^{1/(2+\alpha)}$.
\end{theorem}

The proof of Theorem \ref{thm:bound2} is in Appendix
\ref{apd:margin}.

\begin{Remarks*}

\begin{longlist}[(1)]
\item[(1)] We set the second maximum in (\ref{eqn:noise}) to $-\infty$ if
for an $x$,
$T_0(x,a)$ is constant in $a$ and thus the set
$\mathcal{A}\setminus\arg\max_{a\in\mathcal{A}}T_0(x,a)=\varnothing$.
\item[(2)] %Note that if condition
%(\ref{eqn:noise}) holds for any decomposition of $Q_0(X,A)$ into
%$W_0(X)+T_0(X,A)$, it holds for all.
Condition (\ref{eqn:noise}) is similar to the margin condition in
classification
\cite{polonik1995,mammen1999,tsybakov2004}; in classification this
assumption is often used to obtain
sharp upper bounds on the excess $0$--$1$ risk in terms of other
surrogate risks \cite{bartlett2006}. Here
$\max_{a\in\mathcal{A}}T_0(x, a)-
\max_{a\in\mathcal{A}\setminus\arg\max_{a\in\mathcal{A}}T_0(x,a)}T_0(x,a)$
can be viewed as the ``margin'' of $T_0$ at observation $X=x$. It
measures the difference in mean responses between the optimal
treatment(s) and the best suboptimal treatment(s) at $x$. For example,
suppose $X\sim U[-1,1]$, $P(A=1|X)=P(A=-1|X)=1/2$ and $T_0(X,A)=XA$. Then
the margin condition holds with $C=1/2$ and $\alpha=1$.
Note the margin condition does not exclude multiple optimal treatments
for any observation~$x$.
%This is different from the noise condition in classification. More
% specifically, the goal of classification is to predict the label of
%an observation correctly.
% If $P(Y=1|X=x)=1/2$, it is impossible to correctly tell the class
%label of $X=x$. In the
% decision making problem, the goal is to maximize the mean response.
% If two treatments are equally optimal (i.e. both result in the
%largest mean response) at $x$, we
% do not need to differentiate them.
However, when $\alpha>0$, it does exclude suboptimal treatments that
yield a conditional mean response
very close to the
largest conditional mean response for a set of $x$ with nonzero
probability.

%terms in $Q_0$ that involve $A$. The consequence is that the quality
%of an estimated
%individualized treatment rule only depends on how well we estimate
%$T_0$. Under appropriate conditions, the estimation of $T_0$
%will not be effected by the estimation of $W_0$. See Section

\item[(3)] For $C=1, \alpha=0$, condition (\ref{eqn:noise}) always holds for all
$\epsilon>0$; in this case (\ref{eqn:bound2}) reduces to (\ref{eqn:bound1}).

\item[(4)] The larger the $\alpha$, the larger the exponent $(1+\alpha
)/(2+\alpha)$
and thus the stronger the upper bounds in (\ref{eqn:bound2}) and (\ref
{eqn:bound3}).
However, the
margin condition is unlikely to hold for all $\epsilon$ if $\alpha$ is
very large.
An alternate margin condition and upper bound are as follows.

\textit{Suppose $p(a|x)\geq S^{-1}$ for all $(x,a)$ pairs. Assume there
is an $\epsilon\!>\!0$, such~that
%
%e3.6 ###
\begin{equation}\label{eqn:noise2}
\mathbf{P}\Bigl(\max_{a\in\mathcal{A}}T_0(X,a)
-\max_{a\in\mathcal{A}\setminus\mathop{\arg\max}_{a\in\mathcal{A}}T_0(X,a)}T_0(X,a)
<\epsilon\Bigr)=0.
\end{equation}
Then $V(d_0)-V(d)\leq4S[L(Q)-L(Q_0)]/\epsilon$ and $V(d_0)-V(d)\leq
4SE(T-T_0)^2/\epsilon$.}

The proof is essentially the same as that of Theorem
\ref{thm:bound2} and is omitted. Condition
(\ref{eqn:noise2}) means that $T_0$ evaluated at the optimal
treatment(s) minus $T_0$ evaluated at the best suboptimal treatment(s)
is bounded below by a~positive constant for almost all $X$ observations.
If $X$ assumes only a finite number of values, then this condition always
holds, because we can take $\epsilon$ to be the smallest difference in $T_0$
when evaluated at
the optimal treatment(s) and the suboptimal treatment(s)
[note that if $T_0(x,a)$ is constant for all $a\in\mathcal{A}$ for
some observation $X=x$, then all treatments are optimal for that
observation].

\item[(5)] Inequality (\ref{eqn:bound3}) cannot be improved in the sense
that choosing $T=T_0$ yields zero on both sides of the inequality.
Moreover, an inequality in the opposite direction is not possible,
since each ITR is associated with many
nontrivial $T$-functions. For example, suppose $X\sim U[-1,1]$,
$P(A=1|X)=P(A=-1|X)=1/2$ and $T_0(X,A) = (X-1/3)^2A$. The optimal ITR
is $d_0(X)=1$ a.s. Consider $T(X,A)=\theta A$. Then
maximizing $T(X,A)$ yields the optimal ITR as long as
$\theta>0$. This means that the left-hand side (LHS) of (\ref
{eqn:bound3}) is zero,
while the right-hand side (RHS) is always positive no matter what value
$\theta$
takes.
\end{longlist}
\end{Remarks*}

%{\color{blue} On one hand, Theorem \ref{thm:bound2} supports the
%approach of minimizing
%the empirical quadratic loss to estimate $T_0$ or $Q_0$ and
%then maximizing this estimator over $a\in\mathcal{A}$ to obtain an
%ITR. On the other hand,
%this theorem also indicates a mismatch between the loss functions
%(weighted 0--1 loss and the quadratic loss) because the minimizer of
%the RHS of (\ref{eqn:bound3}) or (\ref{eqn:bound2}) over a class of
%T-, or respectively, Q-functions may not give the ITR that
%minimizes the LHS over the corresponding class of treatment rules.}
%{\color{red}
Theorem \ref{thm:bound2} supports the approach of minimizing
the estimated prediction error to estimate $Q_0$ or $T_0$ and
then maximizing this estimator over $a\in\mathcal{A}$ to obtain an ITR.
It is natural to expect that even when the approximation space
used in estimating $Q_0$ or $T_0$ does not contain the truth, this
approach will
provide the best (highest Value) of the considered ITRs. Unfortunately,
this does not occur due to the mismatch between the loss functions
(weighted 0--1 loss and the quadratic loss). This mismatch is indicated
by remark (5) above.
%It is natural to expect that even when the approximation space
% used in estimating $T_0$ does not contain $T_0$, this approach will
% provide the best (highest value) of the considered ITRs. However,
%Unfortunately
% this does not occur due to the mismatch between the loss functions
% (weighted 0--1 loss and the quadratic loss).
%This is because the minimizer of the right hand side of (
%T-, or respectively, Q-functions may not give the ITR that
%minimizes the left hand side over the corresponding class of treatment
%rules.
% To see this mismatch To see this,
More precisely, note that the approximation space, say $\mathcal{Q}$
for $Q_0$, places implicit restrictions on the
class of ITRs that will be considered. In
effect, the class of ITRs is
$\mathcal{D}_{\mathcal{Q}}=\{d(X)\in\arg\max_aQ(X,a)\dvtx Q\in\mathcal{Q}\}$.
It
turns out that minimizing the prediction error may not result in the
ITR in $\mathcal{D}_{\mathcal{Q}}$ that maximizes the Value. This
occurs when the approximation space $\mathcal{Q}$ does not provide a
treatment effect term close to the treatment effect term in $Q_0$. In
the following toy example, the optimal ITR $d_0$ belongs to
$\mathcal{D}_\mathcal{Q}$, yet the prediction error minimizer over
$\mathcal{Q}$ does not yield $d_0$.
%This mismatch occurs because our
%goal is to maximize the weighted indicator function
%(\ref{eqn:value}) while we minimize the quadratic loss instead.
%
\begin{exam*}
Suppose $X$ is uniformly distributed in $[-1,1]$, $A$ is binary
$\{-1,1\}$ with probability $1/2$ each and is independent of $X$,
and $R$ is normally distributed with mean $Q_0(X,A)=(X-1/3)^2A$ and
variance $1$. It is easy to see that the optimal ITR satisfies
$d_0(X)=1$ a.s. and
$V(d_0)=4/9$.
Consider approximation space $\mathcal{Q}=\{Q(X,A;\bolds{\theta})=(1,
X, A, XA)\bolds{\theta}\dvtx\bolds{\theta}\in\mathbb{R}^4\}$
for $Q_0$. Thus the space of ITRs under consideration is
$\mathcal{D}_{\mathcal{Q}}=\{d(X)=\operatorname{sign}(\theta_3+\theta_4X)\dvtx\theta
_3,\theta_4\in
\mathbb{R}\}$. Note that $d_0\in\mathcal{D}_{\mathcal{Q}}$ since
$d_0(X)$ can be written as $\operatorname{sign}(\theta_3+\theta_4X)$ for any
$\theta_3>0$ and $\theta_4=0$. $d_0$ is the best treatment rule in
$\mathcal{D}_{\mathcal{Q}}$. However, minimizing the prediction
error $L(Q)$ over $\mathcal{Q}$ yields $Q^*(X,A)=(4/9-2/3X)A$. The ITR
associated with $Q^*$ is
$d^*(X)=\arg\max_{a\in\{-1,1\}}Q^*(X,a)=\operatorname{sign}(2/3-X)$, which has lower
Value than $d_0$
($V(d^*)=E[\frac{1_{A(2/3-X)>0}R}{1/2}]=29/81<V(d_0)$).
\end{exam*}

%s4 ###
\section{Estimation via $l_1$-penalized least squares}
\label{sec:lasso}

%From the
% previous section, we see that if the treatment effect term in an
% estimated $\hat Q_0$ is a high quality estimator of $T_0$
%(i.e. the treatment effect term in $Q_0$),
% then the individualized treatment rule, $\hat d_n(X)\in\arg\max_{a\in
%seen, if the treatment $A$ is appropriated coded, then the
%misspecification of the main effect term $E[Q_0(X,A)|X]$ will
%not cause bias in the specification of $T_0$.
%Thus we focus on the estimation of $Q_0$.

To deal with the mismatch between minimizing the
prediction error and maximizing the Value discussed in the prior
section, we consider a large linear approximation space $\mathcal{Q}$
for $Q_0$. Since overfitting is likely (due to the
potentially large number of pretreatment variables and/or large
approximation space for $Q_0$), we use penalized least
squares (see Section S.1 of the supplemental article \cite{supplement}
for further discussion of the overfitting problem).
Furthermore, we use $l_1$-penalized least squares ($l_1$-PLS, \cite
{tibshirani1996}) as the $l_1$ penalty does some
variable selection and as a result will lead to
ITRs that are cheaper to implement (fewer
variables to collect per patient) and easier to interpret. See
Section \ref{sec:discussion} for the discussion of other potential
penalization methods.

Let $\{(X_i,A_i,R_i)\}_{i=1}^n$ represent i.i.d.
observations on $n$ subjects in a randomized trial. For convenience, we use
$E_n$ to denote the associated empirical
expectation [i.e., $E_nf=\sum_{i=1}^nf(X_i,A_i,R_i)/n$ for any
real-valued function
$f$ on $\mathcal{X}\times\mathcal{A}\times\mathbb{R}$].
Let $\mathcal{Q}\triangleq\{Q(X,A;\bolds{\theta})=\Phi(X,A)
\bolds{\theta},\bolds{\theta}\in\mathbb{R}^{J}\}$ be the
approximation space for $Q_0$, where
$\Phi(X,A)=(\phi_1(X,A),\ldots,\phi_{{J}}(X,A))$ is a $1$ by $J$
vector composed of basis functions on
$\mathcal{X}\times\mathcal{A}$, $\bolds{\theta}$ is a $J$ by
$1$ parameter vector, and ${J}$ is the number of basis functions
(for clarity here $J$ will be fixed in $n$, see Appendix
\ref{apd:peoracle} for results with $J$ increasing as $n$
increases).
%For example, $\Phi$ can be a sequence of patient pretreatment
%variables, treatments and their interactions or a sequence of
%wavelet basis functions as described in section
The $l_1$-PLS estimator of
$\bolds{\theta}$ is
%
%e4.1 ###
\begin{equation}
\label{eqn:thetahat}
\bolds{\hat\theta}_n=\mathop{\arg\min}_{\bolds{\theta}\in\mathbb{R}^{J}}
\Biggl\{E_n[R-\Phi(X,A)\bolds{\theta}]^2+{\lambda_n}\sum
_{j=1}^{J}\hat\sigma_j|\theta_j|\Biggr\},
\end{equation}
where $\hat\sigma_j=[E_n\phi_j(X,A)^2]^{1/2}$, $\theta_j$
is the $j$th component of $\bolds{\theta}$ and ${\lambda_n}$
is a tuning parameter that controls the amount of penalization. The
weights $\hat\sigma_j$'s are used to balance the scale of different
basis functions; these weights were used in Bunea, Tsybakov and Wegkamp
\cite{bunea2007} and van de Geer \cite{vandegeer2008}. In some
situations, it is natural to penalize only a subset of coefficients
and/or use different weights in the penalty; see Section S.2 of the supplemental
artic\-le~\cite{supplement} for required modifications. The resulting estimated
ITR satisfies
%
%e4.2 ###
\begin{equation}
\label{eqn:pihat}
\hat d_n(X)\in\mathop{\arg\max}_{a\in\mathcal{A}}\Phi(X,a)\bolds{\hat\theta}_n.
\end{equation}

%s4.1 ###
\subsection{Performance guarantee for the $l_1$-PLS}\label{sec:finaloracle}

In this section, we provide finite sample upper bounds on the
difference between the optimal Value and the Value obtained by the
$l_1$-PLS estimator in terms of the prediction errors resulting from
the estimation of $Q_0$ and
$T_0$. These upper bounds guarantee that if $Q_0$ (or $T_0$) is
consistently estimated, the Value of the estimated ITR will converge
to the optimal Value. Perhaps more importantly,
the finite sample upper bounds provided below do \textit{not}
require the assumption that either $Q_0$ or $T_0$ is consistently
estimated. Thus, each upper bound includes approximation error as well
as estimation error. The estimation error
decreases with decreasing model sparsity and increasing sample size.
An ``oracle'' model for $Q_0$ (or $T_0$) minimizes the sum of these
two errors among suitably sparse linear models [see remark (2) after
Theorem \ref{thm:peoraclefix} for a~precise definition of the oracle
model]. In finite samples, the upper bounds imply that~$\hat d_n$, the ITR
produced by the $l_1$-PLS method, will have Value roughly as if the
$l_1$-PLS method
detects the sparsity of the oracle model and then estimates from the
oracle model using ordinary least squares [see remark (3) below].
%Furthermore in this setting the upper bound can be used to inform how
%to choose the tuning parameter involved in the $l_1$-penalty to
%achieve the best rate of convergence asymptotically.

Define the prediction error minimizer $\bolds{\theta}^*\in\mathbb
{R}^{J}$ by
%
%e4.3 ###
\begin{equation}\label{eqn:thetastarfix}
\bolds{\theta}^*\in\mathop{\arg\min}_{\bolds{\theta}\in\mathbb{R}^{J}}
L(\Phi\bolds{\theta})=\mathop{\arg\min}_{\bolds{\theta}\in\mathbb
{R}^{J}}E(R-\Phi\bolds{\theta})^2.
\end{equation}
For expositional simplicity assume that $\bolds{\theta}^*$ is
unique, and define the sparsity of $\bolds{\theta}\in\mathbb
{R}^J$ by its $l_0$ norm, $\|\bolds{\theta}\|_0$
(see Appendix \ref{apd:peoracle} for a more general setting, where
$\bolds{\theta}^*$ is not unique and a laxer definition of
sparsity is used).
% In that case, we use $[\bolds{\theta}^*]$ to denote the set of
%the most sparse prediction error minimizers, i.e.
%$[\bolds{\theta}^*] = \arg\min_{\bolds{\theta}\in\arg\min_{
% L(\Phi\bolds{\theta})}\|\bolds{\theta}\|_0$, where
%$||\bolds{\theta}||_0$ is the $l_0$ norm of $\bolds{\theta}$.
%{\color{red}
As discussed above, for finite~$n$, instead of estimating $\bolds
{\theta}^*$, the $l_1$-PLS
estimator $\bolds{\hat\theta}_n$ estimates a parameter
$\bolds{\theta}_n^{**}$,
possessing small prediction error and with controlled sparsity.
% }{\color{blue}
%As discussed above, for finite $n$, the $l_1$-PLS
%estimator $\bolds{\hat\theta}_n$ behaves as if it knew the
%sparsity of an oracle model. This oracle model
%is indicated by a parameter, $\bolds{\theta}^{**}_n$, processing
%small prediction error and with controlled sparsity (the oracle model
%contains and only contains basis functions for which the corresponding
%components in $\bolds{\theta}^{**}_n$ are nonzero).}
% {\color{blue} we need to be real explicit about what we mean by
%sparsity--when you deleted the definition of sparsity you made the
%definition of sparsity implicit as opposed to explicit--this causes
%problems with reviewers}
For any bounded function $f$ on $\mathcal{X}\times\mathcal{A}$, let
$\|f\|_\infty
\triangleq\sup_{x\in\mathcal{X},a\in\mathcal{A}}|f(x,a)|$.
$\bolds{\theta}_n^{**}$ lies in the set of parameters $\Theta_n$
defined by\looseness=1
%
%e4.4 ###
\begin{eqnarray} \label{eqn:oracleset}
&&\Theta_{n} \triangleq
\Biggl\{\bolds{\theta}\in\mathbb{R}^{J}\dvtx
\|\Phi(\bolds{\theta}-\bolds{\theta}^*)\|_\infty\leq\eta,\nonumber\\
&&\hspace*{33.7pt}
\max_{j=1,\ldots,J}\biggl|\frac{E[\phi_j\Phi(\bolds{\theta}-
\bolds{\theta}^*)]}{\sigma_j}\biggr|\leq
11\eta\sqrt{\frac{\log(Jn)
}{n}}\\
&&\hspace*{99pt}\mbox{and }\Vert\bolds{\theta}\Vert_0\leq\frac{\beta}{489U}\sqrt{\frac
{n}{\log(Jn)}}
%N_{M_0(\bolds{\theta})}\leq\frac{\beta}{480}
%(\sqrt{\frac{1}{9}+\frac{n}{6U^2\max\{\log2J,\log
%2n\}}}-\frac{1}{3}
%)
\Biggr\},\nonumber
\end{eqnarray}
%
%{\color{red} $$\mbox{Is it better to replace the above condition on
%$N_{M_0(\bolds{\theta})}$ with }
%N_{M_0(\bolds{\theta})}\leq\frac{\beta}{2400U}\sqrt{\frac{n}{\max
%2J,\log2n\}}}\mbox{ ?}$$}
where $\sigma_j=(E\phi_j^2)^{1/2}$, and $\eta$,
$\beta$ and $U$ are positive constants that will be defined in
Theorem \ref{thm:finaloracle}.

The first two conditions in
(\ref{eqn:oracleset}) restrict $\Theta_{n}$ to $\bolds{\theta}$'s
with controlled distance in sup norm and with controlled distance in
%In the first condition in
%(\ref{eqn:oracleset}), $\eta_1$ represents the maximal difference
%between
%$\bolds{\theta}\in\Theta_{n}$ and $\bolds{\theta}^*$ in
%sup norm. In the second condition, $\lambda_n$ represents
%the maximal size of the
prediction error
%of $\bolds{\theta}\in\Theta_{n}$
via first
order derivatives (note that
$|E[\phi_j\Phi(\bolds{\theta}-\bolds{\theta}^*)]/\sigma
_j|
=|\partial L(\Phi\bolds{\theta})/\partial\theta_j-\partial
L(\Phi\bolds{\theta}^*)/\partial\theta_{j}^*|/\break{2\sigma_j}$).
The third condition restricts $\Theta_{n}$ to sparse $\bolds{\theta
}$'s. Note that
as $n$ increases this sparsity requirement becomes laxer,
ensuring that $\bolds{\theta}^*\in\Theta_n$ for sufficiently
large $n$.

When $\Theta_n$ is nonempty, $\bolds{\theta}^{**}_n$ is given by
%
%e4.5 ###
\begin{equation} \label{eqn:thetastarstar}
\bolds{\theta}^{**}_n=\mathop{\arg\min}_{\bolds{\theta}\in\Theta
_n}[L(\Phi\bolds{\theta})
+3\|\bolds{\theta}\|_0\lambda^2_n/\beta].
\end{equation}
Note that $\bolds{\theta}^{**}_n$ is at least
as sparse as $\bolds{\theta}^*$ since by (\ref{eqn:thetastarfix}),
$L(\Phi\bolds{\theta})+3\|\bolds{\theta}\|_0\lambda^2_n/\beta>
L(\Phi\bolds{\theta}^*)+3\|\bolds{\theta}^*\|_0\lambda
^2_n/\beta$ for
any $\bolds{\theta}$ such that
$\|\bolds{\theta}\|_0>\|\bolds{\theta}^*\|_0$.

The following theorem provides a finite sample performance guarantee
for the
ITR produced by the $l_1$-PLS method. Intuitively, this result implies
that if
$Q_0$ can be well approximated by the sparse linear representation
$\bolds{\theta}_n^{**}$ [so
that both $L(\Phi\bolds{\theta}^{**}_n)-L(Q_0)$ and
$\Vert\bolds{\theta}^{**}_n\Vert_0$ are small], then $\hat d_n$ will
have Value
close to the optimal Value in finite samples.
\begin{theorem} \label{thm:finaloracle}
Suppose $p(a|x)\geq S^{-1}$ for a positive constant $S$ for all
$(x,a)$ pairs and the margin condition
(\ref{eqn:noise}) holds for some $C>0$, $\alpha\geq0$ and
all positive~$\epsilon$. Assume:
\begin{longlist}[(1)]
\item[(1)] \hypertarget{ap:errorterm} %prie aplinkos
the error terms $\varepsilon_i=R_i-Q_0(X_i,A_i), i=1,\ldots,n$,
are independent of $(X_i,A_i), i=1,\ldots, n$ and are
i.i.d. with $E(\varepsilon_i)=0$ and
$E[|\varepsilon_i|^l]\leq l!c^{l-2}\sigma^2/2$ for some
$c,\sigma^2>0$ for all $l\geq2$;

\item[(2)] \hypertarget{ap:basis}
there exist finite, positive constants $U$ and $\eta$ such that
$\max_{j=1,\ldots,{J}}$ $\|\phi_j\|_\infty/\sigma_j\leq U$ and
$\|Q_0-\Phi\bolds{\theta}^*\|_\infty\leq\eta$; and

\item[(3)] \hypertarget{ap:grammatrix1} $E[(\phi_1/\sigma_1,\ldots,\phi_J/\sigma_J)^T(\phi_1/\sigma
_1,\ldots,\phi_J/\sigma_J)]$ is positive definite, and the smallest
eigenvalue is denoted by $\beta$.
%there exists a constant $\beta>0$ such that for all
%$\bolds{\theta}$ and $\bolds{\theta}^\prime\in\mathbb{R}^J$
%E[\Phi(\bolds{\theta}^\prime-\bolds{\theta})]^2||\mathbf{
%(\sum_{j\in
%M_0(\bolds{\theta})}\sigma_j|\theta_j^\prime-\theta_j|)^2,
%where $M_0(\bolds{\theta})\triangleq\{j=1,\ldots,J:\theta_j\neq0
\end{longlist}
Consider the estimated ITR $\hat d_n$ defined by (\ref{eqn:pihat}) with tuning
parameter
%
%e4.6 ###
\begin{equation}\label{eqn:lambdaconditionfix}
\lambda_n\geq k \sqrt{\frac{\log(Jn)}{n}},
\end{equation}
where $k=82\max\{c,\sigma,\eta\}$.
%Let $\Theta_n^o$ be defined in
%(\ref{eqn:closeset}). Assume \resume{enumerate}
%there exists a constant $\beta>0$ such that, for all
%$\bolds{\theta}\in\Theta_{n}^o\setminus\{\mathbf{0}\}$ and
%$\mathbf{\tilde\theta}\in\{\mathbb{R}^{J}:
%M_0(\bolds{\theta})}\sigma_j|\tilde\theta_j|\leq
%(2\gamma+5)\sum_{j\in
%M_0(\bolds{\theta})}\sigma_j|\tilde\theta_j-\theta_j|/(1-2\gamma)
%E[\Phi(\mathbf{\tilde\theta}-\bolds{\theta})]^2N_{M_0(
%(\sum_{j\in
%M_0(\bolds{\theta})}\sigma_j|\tilde\theta_j-\theta_j|)^2.
Let $\Theta_n$ be the set defined in (\ref{eqn:oracleset}). Then for
any $n\geq24U^2\log(Jn)$ and for which
$\Theta_n$ is nonempty, we have, with probability at least $1-1/n$,
that
%
%e4.7 ###
\begin{equation} \label{eqn:finaloracle}\qquad
V(d_0)-V(\hat d_n)\leq C'\Bigl[
\min_{\bolds{\theta}\in\Theta_n}\bigl(L(\Phi\bolds{\theta})-L(Q_0)
+3\Vert\bolds{\theta}\Vert_0\lambda_n^2/\beta\bigr)
\Bigr]^{({1+\alpha})/({2+\alpha})},
\end{equation}
where $C'=(2^{2+3\alpha}S^{1+\alpha}C)^{1/(2+\alpha)}$.
\end{theorem}

The result follows from inequality (\ref{eqn:bound2}) in Theorem \ref
{thm:bound2} and inequality (\ref{eqn:peoraclefix}) in Theorem
\ref{thm:peoraclefix}. Similar
results in a more general setting can be obtained by combining
(\ref{eqn:bound2}) with inequality (\ref{eqn:peoracle}) in Appendix
\ref{apd:peoracle}.

\begin{Remarks*}
\begin{longlist}[(1)]
\item[(1)] Note that $\bolds{\theta}^{**}_n$ is the minimizer of the
upper bound on the RHS of (\ref{eqn:finaloracle}) and that
$\bolds{\theta}^{**}_n$ is contained\vadjust{\goodbreak} in the set
$\{\bolds{\theta}^{*,(m)}_n\dvtx m\subset\{1,\ldots,J\}\}$. Each
$\bolds{\theta}^{*,(m)}_n$ satisfies
$\bolds{\theta}^{*,(m)}_n=
\arg\min_{\{\bolds{\theta}\in\Theta_n\dvtx\theta_j=0\ \mathrm{for}\ \mathrm{all}
\ j\notin m\}}L(\Phi\bolds{\theta})$; that is, $\bolds{\theta}^{*,(m)}_n$
minimizes the prediction error of the model indexed by the set $m$
(i.e., model $\{\sum_{j\in m}\phi_j\theta_j\dvtx\theta_j\in\mathbb{R}\}$)
(within $\Theta_n$). For each $\bolds{\theta}^{*,(m)}_n$, the
first term in the upper bound in (\ref{eqn:finaloracle}) [i.e.,
$L(\Phi\bolds{\theta}^{*,(m)}_n)-L(Q_0)$] is the approximation
error of the model indexed by $m$ within~$\Theta_n$. As in\vspace*{1pt}
van de Geer \cite{vandegeer2008}, we call the second term
$3\Vert\bolds{\theta}^{*,(m)}_n\Vert_0\lambda_n^2/\beta$ the estimation
error of the model indexed by $m$. To see why, first put $\lambda_n= k
\sqrt{\log(Jn)/n}$. Then, ignoring the $\log(n)$ factor, the second
term is a function of the sparsity of model $m$ relative to the sample
size, $n$. Up to constants, the second term is a ``tight'' upper bound
for the estimation error of the OLS estimator from model $m$, where
``tight'' means that the convergence rate in the bound is the best
known rate. Note that $\bolds{\theta}_n^{**}$ is the parameter
that minimizes the sum of the two errors over all models. Such a model
(the model corresponding to $\bolds{\theta}_n^{**}$) is called an
oracle model.
%The $\log(n)$ factor in the estimation error is the price paid for not
%knowing the sparsity of the oracle model.
%By using the $l_1$-PLS method, we pay by a factor of $\log(n)$ in the
%estimation error and as an exchange, the $l_1$-PLS estimator behaves
%roughly as if it knew the sparsity of the oracle model and as if is
%was estimated from the oracle model using OLS.
The $\log(n)$ factor in the estimation error can be viewed as the price
paid for not knowing the sparsity of the oracle model and thus having
to conduct model selection.
See remark (2) after Theorem \ref{thm:peoraclefix} for the precise
definition of the oracle model and its relationship to $\bolds
{\theta}_n^{**}$.

\item[(2)] Suppose $\lambda_n = o(1)$. Then in large samples the estimation
error term
$3\Vert\bolds{\theta}\Vert_0\lambda_n^2/\beta$ is negligible. In this
case, $\bolds{\theta}^{**}_n$ is close to
$\bolds{\theta}^*$.
When the model
$\Phi\bolds{\theta}^*$ approximates $Q_0$ sufficiently well, we
see that setting $\lambda_n$ equal to its lower bound in (\ref
{eqn:lambdaconditionfix}) provides the fastest rate of convergence of
the upper bound to zero. More precisely, suppose $Q_0
=\Phi\bolds{\theta}^*$ [i.e.,
$L(\Phi\bolds{\theta}^*)-L(Q_0)=0$]. Then inequality
(\ref{eqn:finaloracle}) implies that $V(d_0)-V(\hat d_n)\leq
O_p( (\log n/n)^{(1+\alpha)/(2+\alpha)})$. A~convergence in mean result is
presented in Corollary \ref{cor:convergence}.

\item[(3)] In finite samples, the estimation error
$3\Vert\bolds{\theta}\Vert_0\lambda_n^2/\beta$ is nonnegligible.
The argument of the minimum in the upper bound (\ref{eqn:finaloracle}),
$\bolds{\theta}^{**}_n$, minimizes prediction error among
parameters with controlled sparsity.
In remark (2) after Theorem \ref{thm:peoraclefix}, we discuss how this
upper bound can be viewed as a tight upper bound for the prediction error of the OLS
estimator from an oracle model in the step-wise model selection setting.
%{\color{blue} is closely related to the submodel with
%balanced approximation capability and sparsity (so that the OLS
%estimator from this submodel has the best prediction performance).
%This submodel is also known as the ``oracle'' model.}
In this sense,
inequality (\ref{eqn:finaloracle}) implies that the treatment rule
produced by the $l_1$-PLS method will have a reduction in Value roughly
as if it
knew the sparsity of the oracle model and were estimated from the
oracle model using OLS.

\item[(4)] Assumptions \hyperlink{ap:errorterm}{(1)}--\hyperlink{ap:grammatrix1}{(3)} in Theorem
\ref{thm:finaloracle} are employed to derive the finite sample
prediction error bound for the $l_1$-PLS estimator
$\bolds{\hat\theta}_n$ defined in (\ref{eqn:thetahat}). Below
we briefly discuss these assumptions.

Assumption \hyperlink{ap:errorterm}{(1)} implicitly implies that the error
terms do not have heavy tails. This condition is often assumed to
show that the sample mean of a~variable is concentrated around its
true mean with a high probability. It is easy to verify that this
assumption holds if each $\varepsilon_i$ is bounded. Moreover, it
also holds for some commonly used error distributions that have
unbounded support, such as the normal or double exponential.

Assumption \hyperlink{ap:basis}{(2)} is also used to show the concentration of
the sample mean around the true mean. It is possible to replace the
boundedness condition by a moment condition similar to assumption
\hyperlink{ap:errorterm}{(1)}. This assumption requires that all basis
functions and the difference between $Q_0$ and its best linear
approximation are bounded. Note that we do not assume $\mathcal{Q}$
to be a good approximation space for $Q_0$. However, if
$\Phi\bolds{\theta}^*$ approximates $Q_0$ well, $\eta$ will
be small, which will result in a smaller upper bound in
(\ref{eqn:finaloracle}). In fact, in the generalized result (Theorem
\ref{thm:peoracle}) we allow $U$ and $\eta$ to
increase in $n$.

Assumption \hyperlink{ap:grammatrix1}{(3)} is employed to avoid collinearity.
In fact, we only need
%
%e4.8 ###
\begin{equation}\label{eqn:grammatrix1}
E[\Phi(\bolds{\theta}^\prime-\bolds{\theta})]^2\Vert\bolds
{\theta}\Vert_0\geq\beta
\biggl(\sum_{j\in
M_0(\bolds{\theta})}\sigma_j|\theta_j^\prime-\theta_j|\biggr)^2
\end{equation}
for
$\bolds{\theta}$, $\bolds{\theta}^\prime$ belonging to a
subset of $\mathbb{R}^J$ (see Assumption \ref{apn:grammatrix}),
where $M_0(\bolds{\theta})\triangleq\{j=1,\ldots,J\dvtx\theta_j\neq0\}$.
Condition (\ref{eqn:grammatrix1}) has been used in van de Geer \cite
{vandegeer2008}.
This condition is also similar to the restricted
eigenvalue assumption in Bickel, Ritov and Tsybakov \cite{bickel2008}
in which
$E$ is replaced by $E_n$, and a fixed design matrix is considered.
Clearly, assumption \hyperlink{ap:grammatrix1}{(3)} is a sufficient condition for
(\ref{eqn:grammatrix1}). In addition, condition
(\ref{eqn:grammatrix1}) is satisfied if the correlation
$|E\phi_j\phi_k|/(\sigma_j\sigma_k)$ is small for all $k\in
M_0(\bolds{\theta})$, $j\neq k$ and a subset of $\bolds{\theta
}$'s (similar results in a fixed design
setting have been proved in Bickel, Ritov and Tsybakov \cite
{bickel2008}. The
condition on correlation is also known as ``mutual coherence''
condition in Bunea, Tsybakov and Wegkamp \cite{bunea2007}). See Bickel,
Ritov and Tsybakov
\cite{bickel2008} for other sufficient conditions for
(\ref{eqn:grammatrix1}).
\end{longlist}
\end{Remarks*}

The above upper bound for $V(d_0)-V(\hat d_n)$ involves
$L(\Phi\bolds{\theta})-L(Q_0)$, which measures how well the
conditional mean function $Q_0$ is approximated by~$\mathcal{Q}$.
As we have seen in Section \ref{sec:relation}, the
quality of the estimated ITR only depends
on the estimator of the treatment effect term $T_0$. Below we
provide a~strengthened result in the sense that the upper bound
depends only on how well we approximate the treatment effect term.

First, we identify terms in the linear model $\mathcal{Q}$ that
approximate $T_0$ (recall that $T_0(X,A)\triangleq
Q_0(X,A)-E[Q_0(X,A)|X]$). Without loss of generality, we rewrite the
vector of basis functions as
$\Phi(X,A)=(\Phi^{(1)}(X),\Phi^{(2)}(X,A))$, where
$\Phi^{(1)}=(\phi_1(X),\ldots,\phi_{J^{(1)}}(X))$ is composed of
all components in $\Phi$ that do not contain $A$ and
$\Phi^{(2)}=(\phi_{J^{(1)}+1}(X,A),\ldots,\phi_{J}(X,A))$ is
composed of all components in $\Phi$ that contain $A$.
%Since $A$
%takes only finite values and the randomization distribution $p(a|x)$
%is known,
Note that $A$ takes only finite values. When the randomization
distribution $p(a|x)$ does not depend on $x$,
we can code~$A$ so that
$E[\Phi^{(2)}(X,A)^T|X]=\mathbf{0}$ a.s. (see Section~\ref{sec:realdata}
and Appendix \ref{sec:simdesign}, for examples). For any
$\bolds{\theta}=(\theta_1,\ldots,\theta_J)^T\in\mathbb{R}^J$,
denote
$\bolds{\theta}^{(1)}=(\theta_1,\ldots,\theta_{J^{(1)}})^T$ and
$\bolds{\theta}^{(2)}=(\theta_{J^{(1)}+1},\ldots,\theta_{J})^T$.
Then $\Phi^{(1)}\bolds{\theta}^{(1)}$ approximates
$E[Q_0(X,A)|X]$ and $\Phi^{(2)}\bolds{\theta}^{(2)}$
approximates $T_0$.% Define
%$M_0^{(1)}(\bolds{\theta})=\{j=1,\ldots,J^{(1)}:\theta_j\neq
%0\}$ and
%$M_0^{(2)}(\bolds{\theta})=\{j=J^{(1)}+1,\ldots,J:\theta_j\neq
%0\}$.

The following theorem implies that if the treatment effect term
$T_0$ can be well approximated by a sparse representation, then
$\hat d_n$ will have Value close to the optimal Value.
\begin{theorem} \label{cor:finaloracletopt}
Suppose $p(a|x)\geq S^{-1}$ for a positive constant $S$ for all
$(x,a)$ pairs and the margin condition
(\ref{eqn:noise}) holds for some $C>0$, $\alpha\geq0$ and
all positive~$\epsilon$. Assume
$E[\Phi^{(2)}(X,A)^T|X]= \mathbf{0}$ a.s. Suppose assumptions
\hyperlink{ap:errorterm}{(1)}--\hyperlink{ap:grammatrix1}{(3)} in Theorem \ref
{thm:finaloracle} hold.
Let $\hat d_n$ be the estimated ITR with $\lambda_n$
satisfying condition (\ref{eqn:lambdaconditionfix}). Let $\Theta_n$
be the set defined in (\ref{eqn:oracleset}). Then for any
$n\geq24U^2\log(Jn)$ and for which $\Theta_n$ is
nonempty, we have, with probability at least $1-1/n$, that
%
%e4.9 ###
\begin{eqnarray} \label{eqn:finaloracle1}
&&V(d_0)-V(\hat d_n)\nonumber\\[-8pt]\\[-8pt]
&&\qquad\leq C'\Bigl[
\min_{\bolds{\theta}\in\Theta_n}\bigl(E\bigl(\Phi^{(2)}\bolds{\theta
}^{(2)}-T_0\bigr)^2
+5\bigl\Vert\bolds{\theta}^{(2)}\bigr\Vert_0\lambda_n^2/\beta\bigr)
\Bigr]^{({1+\alpha})/({2+\alpha})},\nonumber
\end{eqnarray}
where $C'=(2^{2+3\alpha}S^{1+\alpha}C)^{1/(2+\alpha)}$.
\end{theorem}

The result follows from inequality (\ref{eqn:bound3}) in Theorem \ref
{thm:bound2} and inequality (\ref{eqn:peoraclefix1}) in Theorem
\ref{thm:peoraclefix}.

\begin{Remarks*}
\begin{longlist}[(1)]
\item[(1)] Inequality\vspace*{1pt} (\ref{eqn:finaloracle1}) improves inequality (\ref
{eqn:finaloracle}) in the
sense that it guarantees a small reduction in Value of $\hat d_n$
[i.e., $V(d_0)-V(\hat d_n)$] as long as
the treatment effect term $T_0$ is well approximated by a sparse
linear representation; it does not require a good approximation of the
entire conditional mean function $Q_0$. In many situations $Q_0$
may be very complex, but $T_0$ could be very simple. This means that
$T_0$ is much more likely to be well approximated as compared to
$Q_0$ (indeed, if there is no difference between treatments, then
$T_0\equiv0$).

%

%functions
%so that the mean square error
%$E(\Phi^{(2)}\bolds{\theta}^{(2),*}-T_0)^2$ converges to $0$ as
%$n\rightarrow\infty$, where $T_0$ is the treatment effect term in
%$Q_0$ and $\Phi^{(2)}\bolds{\theta}^{(2),*}$ is the best linear
%approximation to $T_0$. Although our theoretical result does not
%require this condition, if this condition does hold then our result
%implies that $V(\hat d_n)$ converges to the optimal value. We refer
%to Barron et al. \cite{barron1999} for general results on the
%construction of approximation
% spaces that guarantee this condition.

\item[(2)] Inequality (\ref{eqn:finaloracle1})
cannot be improved in the sense that if there is no treatment effect
(i.e., $T_0\equiv0$), then both sides of the inequality are zero.
This result implies that minimizing the penalized empirical
prediction error indeed yields high Value (at least asymptotically) if
$T_0$ can be well approximated.

%$T_0$ poorly, this inequality is not suitable for justifying the
%quality of decision rule produced by the $l_1$-PLS method since the
%RHS of (\ref{eqn:finaloracle1}) will be large. Indeed, the
%performance of the $l_1$-PLS method is not clear in this case due to
%the mismatch between minimizing the quadratic error and maximizing
%the value. The $\lambda_n$ targeted for prediction (i.e. minimizing
%quadratic loss) may not give the decision rule with best value (even
%asymptotically). To possibly avoid the mismatch problem, we propose
%to select $\lambda_n$ by maximizing a value estimator in practical
%implementation of the $l_1$-PLS method. See Section \ref{sec:data}
%for detailed demonstration and simulation results.
\end{longlist}
\end{Remarks*}

The following asymptotic result follows from Theorem
\ref{cor:finaloracletopt}. Note that when
$E[\Phi^{(2)}(X,A)^T |X] = \mathbf{0}$ a.s.,
$L(\Phi\bolds{\theta})-L(Q_0) = E[\Phi^{(1)}\bolds{\theta
}^{(1)}- E(Q_0|X)]^2 + E[\Phi^{(2)}\bolds{\theta}^{(2)}-T_0]^2 $.
Thus, the estimation of the treatment effect term $T_0$ is
asymptotically separated from the estimation of the main effect term
$E(Q_0|X)$.\vspace*{1pt}
In this case, $\Phi^{(2)}\bolds{\theta}^{(2),*}$ is the best
linear approximation of
the treatment effect term~$T_0$, where $\bolds{\theta}^{(2),*}$ is
the vector of components in $\bolds{\theta}^*$ corresponding to
$\Phi^{(2)}$.
\begin{corollary} \label{cor:convergence}
Suppose $p(a|x)\geq S^{-1}$ for a positive constant $S$ for all
$(x,a)$ pairs and the margin condition\vspace*{1pt}
(\ref{eqn:noise}) holds for some $C>0$, $\alpha\geq0$ and
all positive $\epsilon$. Assume
$E[\Phi^{(2)}(X,A)^T |X] = \mathbf{0}$ a.s. In addition,\vspace*{2pt} suppose assumptions
\hyperlink{ap:errorterm}{(1)}--\hyperlink{ap:grammatrix1}{(3)} in Theorem \ref{thm:finaloracle}
hold. Let $\hat d_n$ be the estimated
ITR with tuning parameter $\lambda_n=k_1 \sqrt{\log(Jn)/n}$ for a
constant $k_1\geq82\max\{c,\sigma,\eta\}$. If
$T_0(X,A)=\Phi^{(2)}\bolds{\theta}^{(2),*}$, then
\[
V(d_0)-\mathbf{E}[V(\hat d_n)]= O\bigl((\log
n/n)^{(1+\alpha)/(2+\alpha)}\bigr).
\]
\end{corollary}

This result provides a guarantee on
the convergence rate of $V(\hat d_n)$ to the optimal Value. More
specifically, it means that if $T_0$ is correctly approximated, then
the Value of $\hat d_n$ will converge to the optimal Value in mean
at rate at least as fast as $(\log n/n)^{(1+\alpha)/(2+\alpha)}$
with an appropriate choice of $\lambda_n$.

%s4.2 ###
\subsection{Prediction error bound for the $l_1$-PLS estimator}\label{sec:lassooracle}

In this section, we provide a finite sample upper bound for the
prediction error of the $l_1$-PLS estimator~$\bolds{\hat\theta}_n$.
This result is needed to prove Theorem \ref{thm:finaloracle}.
Furthermore, this result strengthens existing literature on $l_1$-PLS
method in prediction. Finite sample prediction error bounds for the
$l_1$-PLS estimator in the random design setting have been provided in
Bunea, Tsybakov and Wegkamp \cite {bunea2007} for quadratic loss, van
de Geer \cite{vandegeer2008} mainly for Lipschitz loss, and
Koltchinskii~\cite{kol2009} for a variety of loss functions. With
regards quadratic loss, Koltchinskii~\cite{kol2009} requires the
response $Y$ is bounded, while both Bunea, Tsybakov and Wegkamp
\cite{bunea2007} and van de Geer \cite {vandegeer2008} assumed the
existence of a sparse $\bolds{\theta}\in\mathbb{R}^{J}$ such that
$E(\Phi\bolds{\theta}-Q_0)^2$ is upper bounded by a quantity that
decreases to $0$ at a certain rate as $n\rightarrow\infty$ (by
permitting $J$ to increase with $n$ so $\Phi$ depends on $n$ as well).
We improve the results in the sense that we do not make such
assumptions (see Appendix \ref{apd:peoracle} for results when $\Phi$,
$J$ are indexed by $n$ and $J$ increases with $n$).

As in the prior sections, the sparsity of
$\bolds{\theta}$ is measured by its $l_0$ norm, $\|\bolds
{\theta}\|_0$ (see the Appendix \ref{apd:peoracle} for proofs with a laxer
definition of sparsity). Recall that the parameter $\bolds{\theta
}^{**}_n$ defined in (\ref{eqn:thetastarstar}) has small prediction
error and controlled sparsity.
%The $l_1$-PLS estimator
%$\bolds{\hat\theta}_n$ estimates a parameter, $\mathbf{
%prediction error and controlled sparsity.
% When $\Theta_n$ is non-empty, $\bolds{\theta}^{**}_n$ lies in $
%[L(\Phi\bolds{\theta})
%+3\|\bolds{\theta}\|_0\lambda^2_n/\beta].
% Note that $\bolds{\theta}^{**}_n$ is at least
%as sparse as $\bolds{\theta}^*$ since by (\ref{eqn:thetastarfix}),
%$L(\Phi\bolds{\theta})+3\|\bolds{\theta}\|_0\lambda^2_n/
%L(\Phi\bolds{\theta}^*)+3\|\bolds{\theta}^*\|_0\lambda^2_n/
%any $\bolds{\theta}$ such that
%$\|\bolds{\theta}\|_0>\|\bolds{\theta}^*\|_0$.
%In the following theorem, we show that
%$L(\Phi\bolds{\hat\theta}_n)\leq
%L(\Phi\bolds{\theta}^{**}_n)
%+3\|\bolds{\theta}^{**}_n\|_0\lambda^2_n/n$ with high probability.
%For any $\bolds{\theta}\in\mathbb{R}^J$, the model determined by $
%containing and only containing basis functions for which the
%corresponding components in $\bolds{\theta}$ are nonzero. As one
%should see in remark 3 below, $3\|\bolds{\theta}\|_0
%estimator from the model determined by $\bolds{\theta}\in
%Thus intuitively, this upper bound implies that the $l_1$-PLS estimator
%$\bolds{\hat\theta}_n$ will have prediction error roughly as if
%the model determined by $\bolds{\theta}^{**}_n$ were known and
%were estimated using OLS.
%
%Thus
%the ITR produced by the OLS estimator from the model determined by
%$\bolds{\theta}^{**}_n$ could be simpler than the rule
%produced by OLS estimator from the model determined by $\mathbf{
%
\begin{theorem} \label{thm:peoraclefix}
Suppose assumptions \hyperlink{ap:errorterm}{(1)}--\hyperlink{ap:grammatrix1}{(3)} in
Theorem \ref{thm:finaloracle} hold. For any $\eta_1\geq0$, let
$\bolds{\hat\theta}_n$ be the $l_1$-PLS estimator defined by
(\ref{eqn:thetahat}) with tuning parameter $\lambda_n$ satisfying
condition\vspace*{1pt} (\ref{eqn:lambdaconditionfix}). Let $\Theta_n$ be the set
defined in (\ref{eqn:oracleset}). Then for any $n\geq
24U^2\log(Jn)$ and for which $\Theta_n$ is
nonempty, we have, with probability at least $1-1/n$, that
%
%e4.10 ###
\begin{equation}\label{eqn:peoraclefix}\qquad
L(\Phi\bolds{\hat\theta}_n)\leq
\min_{\bolds{\theta}\in\Theta_n}\bigl(L(\Phi\bolds{\theta})
+3\|\bolds{\theta}\|_0\lambda_n^2/\beta\bigr)
=L(\Phi\bolds{\theta}^{**}_n)
+3\|\bolds{\theta}^{**}_n\|_0\lambda_n^2/\beta.
\end{equation}

Furthermore, suppose $E[\Phi^{(2)}(X,A)^T|X]= \mathbf{0}$ a.s. Then
with probability at least
$1-1/n$,
%
%e4.11 ###
\begin{equation}\label{eqn:peoraclefix1}
E\bigl(\Phi^{(2)}\bolds{\hat\theta}{}^{(2)}_n-T_0\bigr)^2\leq
\min_{\bolds{\theta}\in\Theta_n}\bigl(E\bigl(\Phi^{(2)}\bolds{\theta
}^{(2)}-T_0\bigr)^2
+5\bigl\|\bolds{\theta}^{(2)}\bigr\|_0\lambda_n^2/\beta\bigr).
\end{equation}
\end{theorem}

The results follow from Theorem \ref{thm:peoracle} in
Appendix \ref{apd:peoracle} with $\rho=0$, $\gamma=1/8$, $\eta_1=\eta
_2=\eta$, $t=\log2n$ and some simple algebra [notice that
assumption \hyperlink{ap:grammatrix1}{(3)} in Theorem \ref{thm:finaloracle} is a
sufficient condition for Assumptions \ref{apn:grammatrix} and \ref
{apn:grammatrix_Topt}].

\begin{Remarks*}
Inequality (\ref{eqn:peoraclefix1}) provides a finite sample upper
bound on the mean square difference between $T_0$ and its estimator.
This result is used to prove Theorem \ref{cor:finaloracletopt}. The
remarks below discuss how
inequality (\ref{eqn:peoraclefix}) contributes to the
$l_1$-penalization literature in prediction.
\begin{longlist}[(1)]
\item[(1)] The conclusion of Theorem \ref{thm:peoraclefix} holds for all
choices of $\lambda_n$ that satisfy~(\ref{eqn:lambdaconditionfix}).
Suppose $\lambda_n=o(1)$. Then
$L(\Phi\bolds{\theta}^{**}_n)-L(\Phi\bolds{\theta}^{*})\rightarrow
0$ as $n\rightarrow\infty$ (since $\|\bolds{\theta}\|_0$ is
bounded). Inequality (\ref{eqn:peoraclefix}) implies that
$L(\Phi\bolds{\hat\theta}_n)-L(\Phi\bolds{\theta}^{*})\rightarrow
0$ in probability. To achieve the best rate of convergence, equal
sign should be taken in (\ref{eqn:lambdaconditionfix}).

\item[(2)]
Note that $\bolds{\theta}^{**}_n$
minimizes
$L(\Phi\bolds{\theta})-L(Q_0)+3\|\bolds{\theta}\|_0\lambda
_n^2/\beta$. Below we de\-monstrate that
the minimum of
$L(\Phi\bolds{\theta})-L(Q_0)+3\|\bolds{\theta}\|_0\lambda
_n^2/\beta$
can be viewed as the approximation error plus a ``tight'' upper
bound of the estimation error of an ``oracle'' in the stepwise model
selection framework [when ``$=$'' is taken in
(\ref{eqn:lambdaconditionfix})]. Here ``tight'' means the
convergence rate in the bound is the best known rate, and ``oracle''
is defined as follows.

Let $m$ denote a nonempty subset of the index set $\{1,\ldots,J\}$.
Then each~$m$ represents a model which uses a nonempty subset of
$\{\phi_1,\ldots,\phi_{J}\}$ as basis functions (there are $2^J-1$ such
subsets). Define
\[
\bolds{\hat\theta}{}^{(m)}_n=\mathop{\arg\min}_{\{\bolds{\theta}\in
\mathbb{R}^J\dvtx\theta_j=0\
\mathrm{for}\ \mathrm{all} \ j\notin m\}}E_n(R-\Phi\bolds{\theta})^2
\]
and
\[
\bolds{\theta}^{*,(m)}=\mathop{\arg\min}_{\{\bolds{\theta}\in
\mathbb{R}^J\dvtx\theta_j=0\
\mathrm{for}\ \mathrm{all}\  j\notin m\}}L(\Phi\bolds{\theta}).
\]
In this setting, an ideal model selection criterion will pick model
$m^*$ such that $L(\Phi\bolds{\hat\theta}{}^{(m^*)}_n)=\inf_m
L(\Phi\bolds{\hat\theta}{}^{(m)}_n)$. $\bolds{\hat\theta}{}^{(m^*)}_n$ is
referred as an ``oracle'' in Massart~\cite{massart2005}. Note that the
excess prediction error of each $\bolds{\hat\theta}{}^{(m)}_n$ can be
written as
\[
L\bigl(\Phi\bolds{\hat\theta}{}^{(m)}_n\bigr)-L(Q_0)=\bigl[L\bigl(\Phi\bolds{\theta}^{*,(m)}\bigr)-L(Q_0)\bigr]
+\bigl[L\bigl(\Phi\bolds{\hat\theta}{}^{(m)}_n\bigr)-L\bigl(\Phi\bolds{\theta}^{*,(m)}\bigr)\bigr],
\]
where the first term is called the approximation error of model $m$
and the second term is the estimation error. It can be shown that
\cite{bartlett2008} for each model $m$ and $x_m>0$, with
probability at least $1-\exp(-x_m)$,
\[
L\bigl(\Phi\bolds{\hat\theta}{}^{(m)}_n\bigr)-L\bigl(\Phi\bolds{\theta}^{*,(m)}\bigr)\leq
\mbox{constant}\times\biggl(\frac{x_m+|m|\log(n/|m|)}{n}\biggr)
\]
under
appropriate technical conditions, where $|m|$ is the cardinality of
the index set $m$. To our knowledge, this is the best rate known so
far. Taking $x_m=\log n+|m|\log J$ and using the union bound
argument, we have with probability at least $1-O(1/n)$,
%
%e4.12 ###
\begin{eqnarray}\label{eqn:esterror1}
&&L\bigl(\Phi_n\bolds{\hat\theta}{}^{(m^*)}_n\bigr)-L(Q_0)\nonumber\\
&&\qquad=
\min_{m} \bigl(\bigl[L\bigl(\Phi\bolds{\theta}^{*,(m)}\bigr)- L(Q_0)\bigr]
+
L\bigl(\Phi\bolds{\hat\theta}{}^{(m)}_n\bigr)-L\bigl(\Phi\bolds{\theta}^{*,(m)}\bigr)\bigr)\nonumber\\
&&\qquad\leq\min_{m} \biggl(\bigl[L\bigl(\Phi\bolds{\theta}^{*,(m)}\bigr)-
L(Q_0)\bigr] + \mbox{constant}\times\frac{|m|\log(Jn)}{n}\biggr)\nonumber\\
&&\qquad= \min_{\bolds{\theta}} \biggl([L(\Phi\bolds{\theta})-
L(Q_0)] + \mbox{constant}\times\frac{\|\bolds{\theta}\|_0\log
(Jn)}{n}\biggr).
\end{eqnarray}

On the other hand, take $\lambda_n$ so that condition
(\ref{eqn:lambdaconditionfix}) holds with ``$=$''.
Equation (\ref{eqn:peoraclefix}) implies that, with probability at least
$1-1/n$,
\[
L(\Phi\bolds{\hat\theta}_n)-L(Q_0) \leq
\min_{\bolds{\theta}\in\Theta_n}
\biggl([L(\Phi\bolds{\theta})-L(Q_0)]+\mbox{constant}\times
\frac{\|\bolds{\theta}\|_0\log(Jn)}{n}\biggr),%\label{eqn:esterror2}
\]
which is essentially (\ref{eqn:esterror1}) with the constraint of
$\bolds{\theta}\in\Theta_n$. (The ``\textit{constant}'' in the
above inequalities may take different values.) Since $\bolds{\theta}=\bolds{\theta}^{**}_n$ minimizes the approximation error plus
a tight upper bound for the estimation error in the oracle model,
within $\bolds{\theta}\in\Theta_n$, we
refer to $\bolds{\theta}^{**}_n$ as an oracle.

\item[(3)] The result can be used to emphasize that
$l_1$ penalty behaves similarly as the $l_0$ penalty. Note that
$\bolds{\hat\theta}_n$ minimizes the empirical prediction error
$E_n(R-\Phi\bolds{\theta})^2$ plus an $l_1$ penalty, whereas
$\bolds{\theta}^{**}_n$ minimizes the prediction error
$L(\Phi\bolds{\theta})$ plus an $l_0$ penalty. We provide an
intuitive connection between these two quantities. First, note that
$E_n(R-\Phi\bolds{\theta})^2$ estimates
$L(\Phi\bolds{\theta})$ and $\hat\sigma_j$ estimates
$\sigma_j$. We use ``$\approx$'' to denote this relationship. Thus,
%
%e4.13 ###
\begin{eqnarray}\label{eqn:thetahatint}
&&
E_n(R-\Phi\bolds{\theta})^2
+
\lambda_n\sum_{j=1}^J\hat\sigma_j|\theta_j|\\
&&\qquad\approx L(\Phi\bolds{\theta})+
\lambda_n\sum_{j=1}^J\sigma_j|\theta_j|\nonumber\\
&&\qquad\leq L(\Phi\bolds{\theta})+
\lambda_n\sum_{j=1}^J\sigma_j|\hat\theta_{n,j}-\theta_j|
+\lambda_n\sum_{j=1}^J\sigma_j|\hat\theta_{n,j}|,\nonumber%
\end{eqnarray}
where $\hat\theta_{n,j}$ is the $j$th component of
$\bolds{\hat\theta}_n$. In Appendix \ref{apd:peoracle}, we show
that for any
$\bolds{\theta}\in\Theta_{n}$,
$\lambda_n\sum_{j=1}^J\sigma_j|\hat\theta_{n,j}-\theta_j|$ is upper
bounded by $\|\bolds{\theta}\|_0\lambda_n^2/\beta$ up to a~constant
with a high probability. Thus, $\bolds{\hat\theta}_n$
minimizes (\ref{eqn:thetahatint}) and
$\bolds{\theta}^{**}_n$ roughly minimizes an upper bound of
(\ref{eqn:thetahatint}).

\item[(4)] The constants involved in the theorem can be improved; we focused
on readability as opposed to
providing the best constants.
%
%natural to leave $\theta_1$ not penalized. In this case,
%$\hat\sigma_j$ is modified to
%$\hat\sigma_j=[E_n\phi_j^2-(E_n\phi_j)^2]^{1/2}$ and $\sigma_j$ is
%modified to $\sigma_j=[E\phi_j^2-(E\phi_j)^2]^{1/2}$. A result
%similar to the conclusion of Theorem \ref{thm:peoraclefix} can be
%obtained after some modifications of the proof and the involved
%constants. As discussed in \cite{vandegeer2008}, this does not
%bring in new theoretical complications. See Appendix B for further
%discussion.
\end{longlist}
\end{Remarks*}

%s5 ###
\section{A practical implementation and an evaluation}\label{sec:data}
In this section, we develop a practical implementation of the $l_1$-PLS
method, compare this method to two commonly used alternatives and
lastly illustrate the method using the motivating data from the
Nefazodone-CBASP trial \cite{keller2000}.

A realistic implementation of $l_1$-PLS method should use a
data-dependent method to select the tuning parameter, $\lambda_n$.
Since the primary goal is to maximize the Value, we select $\lambda_n$
to maximize a cross validated Value estimator. For any ITR~$d$, it is
easy to verify that $E[(R-V(d))1_{A=d(X)}/p(A|X)]=0$. Thus, an unbiased
estimator of $V(d)$ is
\[
E_n
\bigl[1_{A=d(X)}R/p(A|X)\bigr]/E_n \bigl[1_{A=d(X)}/p(A|X)\bigr]
\]
\cite{murphy2001} [recall that the randomization distribution $p(a|X)$
is known]. We
split the data into $10$ roughly equal-sized parts; then for each
$\lambda_n$ we apply
the $l_1$-PLS based method on each $9$ parts of the data to obtain an
ITR, and estimate the Value of this ITR using
the remaining part;
% (i.e. the average of $1_{A=d(X)}R/p(A|X)$ divided
%by the average of $1_{A=d(X)}/p(A|X)$ over the remaining part);
the $\lambda_n$ that maximizes the average of the
$10$ estimated Values is selected. Since the Value of an ITR is
noncontinuous in the parameters, this usually results in a set of candidate
$\lambda_n$'s achieving maximal Value. In the simulations below, the
resulting $\lambda_n$ is nonunique in around $97\%$ of the data sets.
If necessary, as a second step we reduce the set of $\lambda_n$'s by
including only $\lambda_n$'s leading to the ITR's using the least number
of variables. In the simulations below, this second
criterion effectively reduced the number of candidate $\lambda_n$'s
in around $25\%$ of the data sets, however multiple $\lambda_n$'s still
remained in
around $90\%$ of the data sets. This is not surprising since the Value of
an ITR only depends on the relative magnitudes of
parameters in the ITR. In the third step we select the $\lambda_n$
that minimizes the 10-fold cross
validated prediction error estimator from the remaining candidate
$\lambda_n$'s;
that is, minimization of the empirical
prediction error is used as a final tie breaker.

%In Section \ref{sec:simulation}, we compare the above version of
% $l_1$-PLS with two other methods.
% In Section \ref{sec:realdata}, we use data collected from the
%Nefazodone-CBASP trial \cite{keller2000} to illustrate the
%application of $l_1$-PLS.

%s5.1 ###
\subsection{Simulations} \label{sec:simulation}

A first alternative to $l_1$-PLS is to use ordinary least squares
(OLS). The
estimated ITR is $\hat
d_{\mathrm{OLS}}\in\arg\max_a\Phi(X,a)\bolds{\hat\theta}_{\mathrm{OLS}}$ where
$\bolds{\hat\theta}_{\mathrm{OLS}}$ is the OLS
estimator of $\bolds{\theta}$. A second alternative is called
``prognosis prediction'' \cite{kent2002}. Usually this method employs
multiple data
sets, each of which involves one active treatment. Then
the treatment associated with the best predicted prognosis
is selected. We implement this method by
estimating $E(R|X,A=a)$ via least squares with $l_1$
penalization for each treatment group (each $a\in\mathcal{A}$)
separately. The
tuning parameter involved in each treatment group is selected by
minimizing the $10$-fold cross-validated prediction error estimator.
The resulting ITR satisfies $\hat
d_{\mathrm{PP}}(X)\in\arg\max_{a\in\mathcal{A}}\hat E(R|X,A=a)$ where the
subscript ``PP''
denotes prognosis prediction.

For simplicity, we consider binary $A$. All three methods use the same
number of
data points and the same number of basis functions but use these data
points/basis
functions differently. $l_1$-PLS and OLS use all $J$ basis functions to conduct
estimation with all $n$ data points whereas the prognosis prediction
method splits the data into the two treatment groups and uses $J/2$
basis functions to conduct estimation with the $n/2$ data points in
each of the two treatment groups. To ensure the comparison is fair
across the three methods, the approximation model for each treatment
group is consistent with the approximation model used in both $l_1$-PLS
and OLS [e.g., if $Q_0$ is approximated by $(1, X, A,
XA)\bolds{\theta}$ in $l_1$-PLS and OLS, then in prognosis
prediction we approximate
$E(R|X,A=a)$ by $(1, X)\bolds{\theta}_{\mathrm{PP}}$ for each treatment
group].
We do not penalize the intercept coefficient in either prognosis
prediction or $l_1$-PLS.

The three methods are compared using two criteria: (1)
Value maximization; and (2) simplicity of the estimated ITRs (measured
by the number of variables/basis functions used in
the rule).

We illustrate the comparison of the three methods using $4$ examples
selected to reflect three
scenarios (see Section S.3 of the supplemental article \cite
{supplement} for $4$ further examples):
\begin{longlist}[(1)]
\item[(1)] There is no treatment effect [i.e., $Q_0$ is constructed so that
$T_0=0$; example~(1)]. In this case, all ITRs yield the same Value. Thus, the
simplest rule is preferred.

\item[(2)] There is a treatment effect and the treatment effect term $T_0$
is correctly modeled [example (4) for large $n$ and example (2)]. In this
case, minimizing the prediction error will yield the
ITR that maximizes the Value.

\item[(3)] There is a treatment effect and the treatment effect term
$T_0$ is misspecified [example (4) for small $n$ and example (3)]. In this
case, there might be
a mismatch between prediction error minimization and Value
maximization.
\end{longlist}

The examples are generated as follows. The
treatment $A$ is generated uniformly from $\{-1,1\}$ independent of $X$
and the response $R$. The response $R$ is
normally distributed with mean $Q_0(X,A)$. In examples (1)--(3),
$X\sim U[-1,1]^5$ and we consider three simple examples for $Q_0$. In
example (4),
$X \sim U[0,1]$ and we use a complex $Q_0$, where
$Q_0(X,1)$ and $Q(X, -1)$ are similar to the blocks function used in
Donoho and Johnstone
\cite{donoho1994}.
%To make the simulations more
%realistic, examples (5)-(8) are based on data from the Nefazodone-CBASP
%trial \cite{keller2000} (see Section \ref{sec:realdata} for
%description of the trial). We consider $50$ pretreatment variables
%collected from the trial (i.e. $X\in\mathbb{R}^{50}$) and four
%examples for $Q_0$.
Further details of the simulation design are
provided in Appendix \ref{sec:simdesign}.

\begin{figure}
\vspace{-2pt}
\includegraphics{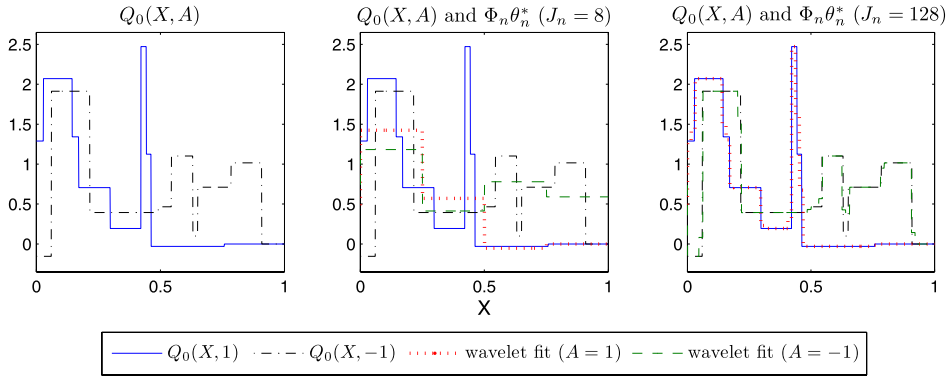}
\vspace{-5pt}
\caption{Plots for:
the conditional mean function $Q_0(X,A)$ (\textup{left}), $Q_0(X,A)$ and the
associated best wavelet fit when $J_n=8$ (\textup{middle}), and $Q_0(X,A)$
and the associated best wavelet fit when $J_n=128$ (\textup{right}) [example
(4)].}
% with parameters specified in (\ref{eqn:parameter}).
\label{fig:qoptfit}
\vspace{-8pt}
\end{figure}

We consider two types of approximation models for $Q_0$. In examples
(1)--(3), we approximate $Q_0$ by
$(1,X,A,XA)\bolds{\theta}$. In example
(4), we approximate $Q_0$ by Haar wavelets. The number of basis
functions may increase as $n$ increases (we index $J$, $\Phi$ and
$\bolds{\theta}^*$ by $n$ in this case). Plots for $Q_0(X,A)$
and the associated best wavelet fits
$\Phi_n(X,A)\bolds{\theta}^*_n$ are provided in Figure
\ref{fig:qoptfit}.

%
%f1 ###

For each example, we simulate data sets of sizes $n=2^k$ for $k=5,\ldots,10$.
%For each of the examples 5 - 8, we simulate
%data sets of size $n=500$.
$1\mbox{,}000$ data sets are generated for each
sample size.
%We apply $l_1$-PLS (denoted by
%$l_1$-PLS), the OLS method (denoted by OLS) and the method based on
%separate prognosis prediction for each treatment (denoted by PP) on
%each data set.
The Value of each estimated ITR is evaluated via Monte Carlo using a
test set of size $10\mbox{,}000$.
The Value of the optimal ITR is also evaluated using the
test set.

%
%f2 ###
\begin{figure}
\vspace{-2pt}
\includegraphics{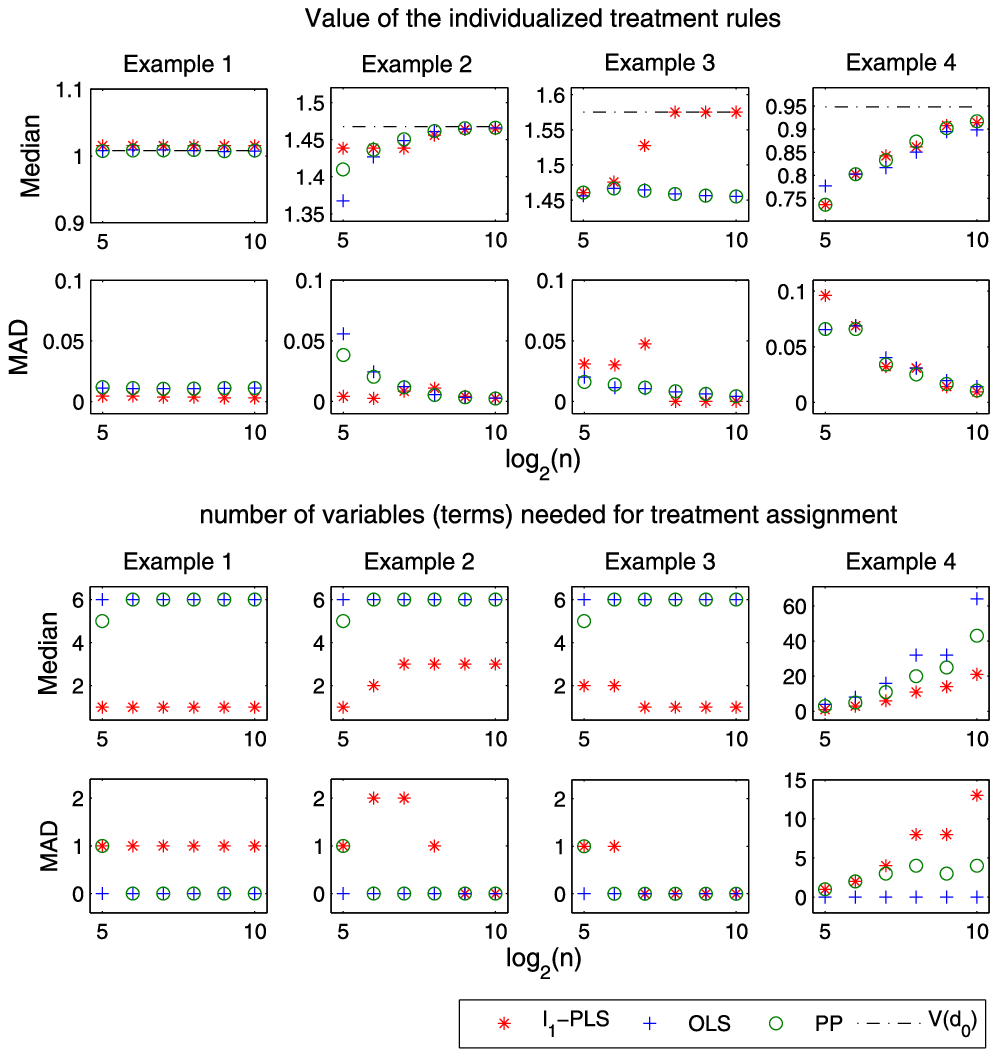}
\vspace{-5pt}
\caption{Comparison of the $l_1$-PLS based method
with the OLS method and the PP method [examples (1)--(4)]: plots for
medians and median absolute deviations (MAD) of
the Value of the estimated decision rules (top panels) and the number
of variables (terms) needed for treatment assignment (including the
main treatment
effect term, bottom panels) over $1\mbox{,}000$ samples versus sample size
on the log scale. The black dash-dotted line in each plot on the
first row denotes the Value of the optimal treatment rule,
for each example. [$n=32, 64, 128, 256, 512, 1024$. The
corresponding numbers of basis functions in example (4) are
$J_n=8, 16, 32, 64, 64, 128$.]}\label{fig:simulation}
\vspace{-8pt}
\end{figure}

Simulation results are presented in Figure \ref{fig:simulation}.
% and
%Table \ref{table:simulation}.
When the approximation model is of
high quality, all methods produce ITRs with similar Value [see examples
(1), (2) and example (4) for large $n$].
However, when the approximation model is poor, the $l_1$-PLS method
may produce highest Value [see example (3)]. Note
that in example (3) settings in which the sample size is small, the Value
of the ITR produced by
$l_1$-PLS method has larger median absolute deviation (MAD)
than the other two methods. One
possible reason is that due to the mismatch between maximizing the
Value and minimizing the prediction error, the Value estimator
plays a strong role in selecting $\lambda_n$.
%{\color{red} I found that the distribution of the value of the
%estimated ITR is bimodal. One mode is at the optimal value and the
%other mode is at the value produced by the prediction error minimizer.}
The nonsmoothness of the Value estimator combined with the mismatch results
in very different $\lambda_n$'s and thus the estimated decision rules
vary greatly
from data set to data set in this example. Nonetheless, the $l_1$-PLS
method is still preferred after taking the variation into account; indeed
$l_1$-PLS produces ITRs with higher Value than both OLS
and PP in around $46\%$, $55\%$ and $67\%$ in data sets of sizes $n=32, 64$
and $128$, respectively. Furthermore, in general the $l_1$-PLS
method uses much fewer variables for treatment assignment than the
other two methods. This is expected because the OLS method
does not have variable selection functionality and the PP method
will use all variables that are predictive of the response $R$ whereas
the use of the Value in selecting the tuning parameter in $l_1$-PLS
discounts variables that are only useful in
predicting the response (and less useful in selecting the best treatment).

\subsection{Nefazodone-CBASP trial example} \label{sec:realdata}

The Nefazodone-CBASP trial was conducted to compare the efficacy of
several alternate treatments for patients with chronic depression.
The study randomized $681$ patients with nonpsychotic chronic major
depressive disorder (MDD) to either Nefazodone, cognitive
behavioral-analysis system of psychotherapy (CBASP) or the
combination of the two treatments. Various assessments were taken
throughout the study, among which the score on the 24-item Hamilton
Rating Scale for Depression (HRSD) was the primary outcome. Low HRSD
scores are desirable. See Keller et al. \cite{keller2000} for more
detail of the
study design and the primary analysis.

In the data analysis, we use a subset of the Nefazodone-CBASP data
consisting of $656$ patients for whom the response HRSD score was
observed. In this trial, pairwise comparisons show that the
combination treatment resulted in significantly lower HRSD scores
than either of the single treatments. There was no overall
difference between the single treatments.
%( than either of the two single treatments will result in much lower
%HRSD
%scores than nefazodone alone or CBASP alone (p-value $<0.001$ for
%both comparisons; the sample mean HRSD score is $10.9505$ in the
%combination treatment group, $15.7569$ in nefazodone group and
%$15.7870$ in CBASP group).

We use $l_1$-PLS to develop an ITR. In
the analysis, the HRSD score is reverse coded so that higher is
better. We consider $50$ pretreatment variables
$X=(X_1,\ldots,X_{50})$. Treatments are coded using contrast coding
of dummy variables $A = (A_1,A_2)$, where $A_1 = 2$ if the
combination treatment is assigned and $-1$ otherwise and $A_2 = 1$
if CBASP is assigned, $-1$ if nefazodone and $0$ otherwise. The
vector of basis functions, $\Phi(X,A)$, is of the form
$(1,X,A_1,XA_1, A_2, XA_2)$. So the number of basis functions is
$J=153$. As a contrast, we also consider the OLS method and the PP
method (separate prognosis prediction for each treatment). The
vector of basis functions used in PP is $(1, X)$ for each of the three treatment
groups. Neither the intercept term nor the main treatment effect
terms in $l_1$-PLS or PP is penalized (see Section S.2 of the
supplemental article \cite{supplement} for the modification of the
weights $\hat\sigma_j$
used in (\ref{eqn:thetahat})).

The ITR given by the $l_1$-PLS method
recommends the combination treatment to all (so none of the
pretreatment variables enter the rule). On the other hand, the PP
method produces an ITR that uses $29$ variables. If the
rule produced by PP were used to assign
treatment for the $656$ patients in the trial, it would recommend
the combination treatment for $614$ patients and nefazodone for the
other $42$ patients. In addition, the OLS method will use all the
$50$ variables. If the ITR produced by OLS were used to
assign treatment for the $656$ patients in the trial, it would
recommend the combination treatment for $429$ patients, nefazodone
for $145$ patients and CBASP for the other $82$ patients.

%We have found that, in general, if one treatment is overwhelmingly
%better than the other treatments, the treatment rules produced by
%both methods are likely to recommend the same treatments for most
%patients; however as the difference in treatments decreases these
%two treatment rules will recommend different treatments for more and
%more patients. To see this, we consider the following 5 examples.
%The first example uses the original data, in which the combination
%treatment is overwhelmingly better. Cohen's f effect size index is
%around $0.25$ (Cohen's f index is the square root of the
%between-group variance divided by the square root of the
%within-group variance; $0.25$ is considered as a medium effect size
%constant from the reverse coded HRSD scores for the combination
%treatment group so that the Cohen's f index is around $0.2$, $0.15$,
%$0.1$ and $0.05$, respectively. Both methods are used on each
%example. The $l_1$-PLS method produces treatment rules that use $0,
%0, 0, 30$ and $23$ variables and the PP method produces treatment
%rules that always use $29$ variables for treatment assignment for
%examples 1 to 5, respectively. If the treatment rules produced by
%the two methods were used to assign treatment for the $656$ patients
%in the trial, they would recommend different treatments on $42$,
%$81$, $132$, $264$ and $331$ patients for examples 1 to 5,
%respectively.

%s6 ###
\section{Discussion}\label{sec:discussion}

Our goal is to construct a high quality ITR that will benefit future patients.
We considered an $l_1$-PLS
based method and provided a finite sample upper bound for
$V(d_0)-V(\hat d_n)$, the reduction in Value of the estimated ITR.

The use of an $l_1$ penalty allows us
to consider a large model for the conditional mean function $Q_0$
yet permits a sparse estimated ITR. In
fact, many other penalization methods such as SCAD \cite{fan2001}
and $l_1$ penalty with adaptive weights (adaptive Lasso;
\cite{zou2006}) also have this property. We choose the nonadaptive
$l_1$ penalty to represent these methods. Interested readers may
justify other PLS methods using similar proof techniques.

%An important issue is how to select the sequence of basis functions
%so that the mean square error
%$E(\Phi^{(2)}\bolds{\theta}^{(2),*}-T_0)^2$ converges to $0$ as
%$n\rightarrow\infty$, where $T_0$ is the treatment effect term in
%$Q_0$ and $\Phi^{(2)}\bolds{\theta}^{(2),*}$ is the best linear
%approximation of $T_0$. Although our theoretical result does not
%require this condition, if this condition does hold then our result
%implies that $V(\hat d_n)$ converges to the optimal value. We refer
%to Barron et al. \cite{barron1999} for general results on the
%construction of approximation
% spaces that guarantee this condition.
The high\vspace*{1pt} probability finite sample upper bounds [i.e., (\ref
{eqn:finaloracle}) and
(\ref{eqn:finaloracle1})] cannot be used to construct a prediction/confidence
interval for $V(d_0)-V(\hat d_n)$ due to the unknown quantities in
the bound. How to develop a tight computable upper bound to
assess the quality of $\hat d_n$ is an open question.

We used cross validation with Value maximization to select the
tuning parameter involved in the $l_1$-PLS method. As compared to
the OLS method and the PP method, this method
may yield higher Value when $T_0$ is misspecified.
However, since only the Value is used to select the tuning
parameter, this method may produce a complex ITR for which the Value is
only slightly higher than that
of a much simpler ITR. In this case, a simpler rule may
be preferred due to the interpretability and cost of collecting the
variables. Investigation of a tuning parameter selection criterion
that trades off the Value with the number of variables in an
ITR is needed.

This paper studied a one stage decision problem. However, it is
evident that some diseases require time-varying treatment. For
example, individuals with a chronic disease often experience a
waxing and waning course of illness. In these settings, the goal is
to construct a sequence of ITRs that
tailor the type and dosage of treatment through time according to an
individual's changing status. There is an abundance of statistical
literature in this area
\cite
{thall2000,thall2002,murphy2003,murphy2005,robins2004,lunceford2002,vanderlaan2005,wahedtsiatis2006}.
Extension of the least squares
based method to the multi-stage decision problem has been presented
in Murphy \cite{murphy2005}. The performance of $l_1$ penalization
in this setting is unclear and worth investigation.

\begin{appendix}
%% redefine the command that creates the equation no.

\section*{Appendix}

%s6.1 ###
\subsection{\texorpdfstring{Proof of Theorem \lowercase{\protect\ref{thm:bound2}}}{Proof of Theorem 3.1}}
\label{apd:margin}

For any ITR $d\dvtx\mathcal{X}\rightarrow\mathcal{A}$, denote
$\triangle
T_d(X)$ $\triangleq\max_{a\in\mathcal{A}}T_0(X,a)-T_0(X,d(X))$. Using
similar arguments to that in Section~\ref{sec:prelim}, we have
$V(d_0)-V(d)=E(\triangle T_d)$. If $V(d_0)-V(d)=0$, then
(\ref{eqn:bound2}) and (\ref{eqn:bound3}) automatically hold.
Otherwise, $E(\triangle T_d)^2\geq(E\triangle T_d)^2>0$. In this
case, for any $\epsilon>0$, define the event\vspace{-2pt}
\[
\Omega_\epsilon=\Bigl\{\max_{a\in\mathcal{A}}T_0(X,a)-
\max_{a\in\mathcal{A}\setminus\mathop{\arg\max}_{a\in\mathcal
{A}}T_0(X,a)}T_0(X,a)\leq\epsilon\Bigr\}.
\]
Then $\triangle T_d\leq(\triangle T_d)^2/\epsilon$ on the event
$\Omega_\epsilon^C$. This together with the fact that $\triangle
T_d\leq(\triangle T_d)^2/\epsilon+\epsilon/4$ implies
\begin{eqnarray*}
V(d_0)-V(d) &=& E(1_{\Omega_\epsilon^C}\triangle T_d)+
E(1_{\Omega_\epsilon}\triangle T_d)\\[-2pt]
&\leq&\frac{1}{\epsilon}E[1_{\Omega_\epsilon^C} (\triangle
T_d)^2]+E\biggl[1_{\Omega_\epsilon}\biggl(\frac{(\triangle
T_d)^2}{\epsilon}+\frac{\epsilon}{4}\biggr)\biggr]\\[-2pt]
&=&\frac{1}{\epsilon}E[(\triangle
T_d)^2]+\frac{\epsilon}{4}P(\Omega_\epsilon)\leq
\frac{1}{\epsilon}E[(\triangle
T_d)^2]+\frac{C}{4}\epsilon^{1+\alpha},
\end{eqnarray*}
where the last inequality follows from the margin condition
(\ref{eqn:noise}). Choosing $\epsilon= (4E(\triangle
T_d)^2/C)^{1/(2+\alpha)}$ to minimize the above upper bound
yields
%
%e6.1 ###
\begin{equation} \label{eqn:thm1part1}
V(d_0)-V(d)\leq
2^{\alpha/(2+\alpha)}C^{1/(2+\alpha)}[E(\triangle
T_d)^2]^{(1+\alpha)/(2+\alpha)}.
\end{equation}

Next, for any $d$ and $Q$ such that
$d(X)\in\max_{a\in\mathcal{A}}Q(X,a)$, let $T(X,A)$ be the associated
treatment effect term. Then
\begin{eqnarray*}
E(\triangle T_d)^2
&=&E\Bigl[\Bigl(\max_{a\in\mathcal{A}}T_0(X,a)-\max_{a\in\mathcal{A}}T(X,a)+
T(X,d(X))-T_0(X,d(X))\Bigr)^2\Bigr]\\[-2pt]
&\leq&
2E\Bigl[\Bigl(\max_{a\in\mathcal{A}}T_0(X,a)-\max_{a\in\mathcal
{A}}T(X,a)\Bigr)^2\\[-2pt]
&&\hspace*{17.8pt}{} +
\bigl(T(X,d(X))-T_0(X,d(X))\bigr)^2\Bigr]\\[-2pt]
&\leq&
4E\Bigl[\max_{a\in\mathcal{A}}\bigl(T(X,a)-T_0(X,a)\bigr)^2\Bigr],
\end{eqnarray*}
where the last inequality follows from the fact that neither
$|{\max_aT_0(X,a)}-\max_aT(X,a)|$ nor $|T(X,d(X))-T_0(X,d(X))|$ is
larger than ${\max_a}|T(X,a)-T_0(X,a)|$. Since $p(a|x)\geq S^{-1}$ for
all $(x,a)$ pairs, we have
%
%e6.2 ###
\begin{eqnarray} \label{eqn:thm1part2}
E(\triangle T_d)^2 &\leq&4S
E\Bigl[\sum_{a\in\mathcal{A}}\bigl(T(X,a)-T_0(X,a)\bigr)^2p(a|X)
\Bigr]\nonumber\\[-2pt]
&=&4SE\bigl(T(X,A)-T_0(X,A)\bigr)^2.
\end{eqnarray}
Inequality (\ref{eqn:bound3}) follows by substituting
(\ref{eqn:thm1part2}) into (\ref{eqn:thm1part1}). Inequality
(\ref{eqn:bound2}) can be proved similarly by noticing that
$\triangle T_d(X)=\max_{a\in\mathcal{A}}Q_0(X,a)-Q_0(X,d(X))$.\vspace{-2pt}

%% redefine the command that creates the equation no.

%s6.2 ###
\subsection{\texorpdfstring{Generalization of Theorem \lowercase{\protect\ref{thm:peoraclefix}}}%
{Generalization of Theorem 4.3}}
\label{apd:peoracle}

In this section, we present a generalization of Theorem
\ref{thm:peoraclefix} where $J$ may depend on $n$ and the sparsity
of any $\bolds{\theta}\in\mathbb{R}^{J}$ is measured by the
number of ``large'' components in $\bolds{\theta}$ as described
in Zhang and Huang \cite{chzhang2008}. In this case, $J$, $\Phi$ and
the prediction
error minimizer $\bolds{\theta}^*$
are denoted as $J_n, \Phi_n$ and $\bolds{\theta}^*_n$,
respectively. All relevant quantities and assumptions are restated below.

Let $|M|$ denote the cardinality of any index set
$M\subseteq\{1,\ldots,J_n\}$. For any $\bolds{\theta}\in
\mathbb{R}^{J_n}$ and constant $\rho\geq0$, define
\[
M_{\rho\lambda_n}(\bolds{\theta})\in\mathop{\arg\min}_{\{M\subseteq\{
1,\ldots,J_n\}\dvtx
\sum_{j\in\{1,\ldots,J_n\}\setminus M}\sigma_j|\theta_j|\leq\rho
|M|\lambda_n\}}|M|.
\]
Then $M_{\rho\lambda_n}(\bolds{\theta})$ is the smallest index
set that contains only ``large'' components
in~$\bolds{\theta}$. $|M_{\rho\lambda_n}(\bolds{\theta})|$
measures the sparsity of $\bolds{\theta}$. It is easy to see that
when $\rho=0$,
$M_0(\bolds{\theta})$ is the index set of nonzero components in
$\bolds{\theta}$ and $|M_0(\bolds{\theta})|=\|\bolds{\theta}\|_0$. Moreover,
$M_{\rho\lambda_n}(\bolds{\theta})$ is an empty set if and only
if $\bolds{\theta}=\mathbf{0}$.

Let $[\bolds{\theta}^*_n]$ be the set of most sparse prediction error
minimizers in the linear model, that is,\vspace{-5pt}
%
%e6.3 ###
\begin{equation}\label{eqn:thetastar}
[\bolds{\theta}^*_n]=\mathop{\arg\min}_{\bolds{\theta}\in\mathop{\arg\min}
_{\bolds{\theta}}
L(\Phi_n\bolds{\theta})}|M_{\rho\lambda_n}(\bolds{\theta})|.
\end{equation}
Note that $[\bolds{\theta}^*_n]$ depends on $\rho\lambda_n$.
% In that case, we use $[\bolds{\theta}^*]$ to denote the set of
%the most sparse prediction error minimizers, i.e.
%$[\bolds{\theta}^*] = \arg\min_{\bolds{\theta}\in\arg\min_{
% L(\Phi\bolds{\theta})}\|\bolds{\theta}\|_0$, where
%$||\bolds{\theta}||_0$ is the $l_0$ norm of $\bolds{\theta}$.

To derive the finite sample upper bound for $L(\Phi_n\bolds{\hat
\theta}_n)$, we need the following assumptions.\vspace{-5pt}
\begin{assumption}\label{apn:errorterm} The error\vspace*{1pt} terms $\varepsilon_i,
i=1,\ldots,n$ are independent of $(X_i,A_i), i=1,\ldots,n$ and are
i.i.d. with $E(\varepsilon_i)=0$ and
$E[|\varepsilon_i|^l]\leq\frac{l!}{2}c^{l-2}\sigma^2$ for some
$c,\sigma^2>0$ for all $l\geq2$.
\end{assumption}
\begin{assumption}\label{apn:basis} For all $n\geq1$:

\begin{longlist}[(a)]
\item[(a)]
there exists an $1\leq U_n<\infty$ such that
$\max_{j=1,\ldots,{J_n}}\|\phi_j\|_\infty/\sigma_j\leq U_n$, where
$\sigma_j\triangleq(E\phi_j^2)^{1/2}$.

\item[(b)] there exists an $0<\eta_{1,n}<\infty$, such that
$\sup_{\bolds{\theta}\in
[\bolds{\theta}^*_n]}\|Q_0-\Phi_n\bolds{\theta}\|_\infty\leq
\eta_{1,n}$.
\end{longlist}
\end{assumption}

For any $0\leq\gamma<1/2$, $\eta_{2,n}\geq0$ (which may
depend on $n$) and tuning parameter~$\lambda_n$, define
\begin{eqnarray*}
\Theta_{n}^o &=&
\biggl\{\bolds{\theta}\in\mathbb{R}^{J_n}\dvtx\exists
\bolds{\theta}^o\in[\bolds{\theta}^*_n]
\mbox{ s.t. }
\|\Phi_n(\bolds{\theta}-\bolds{\theta}^o)\|_\infty\leq\eta
_{2,n}\\
&&\hspace*{31.13pt}
\mbox{ and }
\max_{j=1,\ldots,J_n}\biggl|E\biggl[\Phi_n(\bolds{\theta}-\bolds{\theta}^o)\frac{\phi_j}{\sigma_j}\biggr]\biggr|\leq
\gamma{\lambda_n}\biggr\}.
\end{eqnarray*}
\begin{assumption}\label{apn:grammatrix} For any $n\geq1$,
there exists a $\beta_n>0$ such that
\[
E[\Phi_n(\bolds{\tilde\theta}-\bolds{\theta})]^2|M_{\rho
\lambda_n}(\bolds{\theta})|\geq\beta_n
\biggl[\biggl(\sum_{j\in
M_{\rho\lambda_n}(\bolds{\theta})}\sigma_j|\tilde\theta_j-\theta
_j|\biggr)^2-\rho^2
|M_{\rho\lambda_n}(\bolds{\theta})|^2\lambda_n^2\biggr]
\]
for all
$\bolds{\theta}\in\Theta_{n}^o\setminus\{\mathbf{0}\}$,
$\bolds{\tilde\theta}\in\mathbb{R}^{J_n}$ satisfying $\sum_{j\in
\{1,\ldots,J_n\}\setminus
M_{\rho\lambda_n}(\bolds{\theta})}\sigma_j|\tilde\theta_j|\leq
\frac{2\gamma+5}{1-2\gamma}\times(\sum_{j\in
M_{\rho\lambda_n}(\bolds{\theta})}\sigma_j|\tilde\theta_j-\theta
_j|+\rho
|M_{\rho\lambda_n}(\bolds{\theta})|\lambda_n)$.
\end{assumption}

When $E(\Phi_n^{(2)}(X,A)^T|X)=\mathbf{0}$ a.s. ($\Phi_n^{(2)}$
is defined in Section \ref{sec:finaloracle}), we need an extra
assumption to derive the finite sample upper bound for the mean square
error of the treatment effect estimator $E[\Phi_n^{(2)}\bolds
{\hat\theta}{}^{(2)}_n-T_0(X,A)]^2$ (recall that $T_0(X,A)\triangleq
Q_0(X,A)-E[Q_0(X,A)|X]$).
\begin{assumption}
\label{apn:grammatrix_Topt} For any $n\geq1$, there exists a
$\beta_n>0$ such that
\begin{eqnarray*}
&&E\bigl[\Phi_n^{(2)}\bigl(\bolds{\tilde\theta}{}^{(2)}-\bolds{\theta}^{(2)}\bigr)\bigr]^2
\bigl|M_{\rho\lambda_n}^{(2)}(\bolds{\theta})\bigr|\\
&&\hspace*{0pt}\qquad\geq\beta_n
\biggl[\biggl(\sum_{j\in
M_{\rho\lambda_n}^{(2)}(\bolds{\theta})}\sigma_j|\tilde\theta
_j-\theta_j|\biggr)^2-\rho^2
\bigl|M_{\rho\lambda_n}^{(2)}(\bolds{\theta})\bigr|^2\lambda_n^2\biggr]
\end{eqnarray*}
for all
$\bolds{\theta}\in\Theta_{n}^o\setminus\{\mathbf{0}\}$,
$\bolds{\tilde\theta}\in\mathbb{R}^{J_n}$ satisfying $\sum_{j\in
\{1,\ldots,J_n\}\setminus
M_{\rho\lambda_n}(\bolds{\theta})}\sigma_j|\tilde\theta_j|\leq
\frac{2\gamma+5}{1-2\gamma}\times(\sum_{j\in
M_{\rho\lambda_n}(\bolds{\theta})}|\tilde\theta_j-\theta_j|+\rho
|M_{\rho\lambda_n}(\bolds{\theta})|\lambda_n)$, where
\[
M_{\rho\lambda_n}^{(2)}(\bolds{\theta})\in\mathop{\arg\min}_{\{M\subseteq\{
J_n^{(1)}+1,\ldots,J_n\}\dvtx
\sum_{j\in\{J_n^{(1)}+1,\ldots,J_n\}\setminus
M}\sigma_j|\theta_j|\leq\rho|M|\lambda_n\}}|M|
\]
is the smallest index set that contains only large components in
$\bolds{\theta}^{(2)}$.
\end{assumption}

Without loss of generality, we assume that Assumptions \ref
{apn:grammatrix} and \ref{apn:grammatrix_Topt} hold with the same value
of $\beta_n$. And we can always choose a small enough $\beta_n$ so that
$\rho\beta_n\leq1$ for a given $\rho$.

For any given $t>0$, define
%
%e6.4 ###
\begin{eqnarray}\label{eqn:thetan2}
\Theta_{n}&=&
\Biggl\{\bolds{\theta}\in\Theta_{n}^o\dvtx|M_{\rho\lambda
_n}(\bolds{\theta})|\nonumber\\[-8pt]\\[-8pt]
&&\hspace*{6.5pt}\leq
\frac{(1-2\gamma)^2\beta_n}{120}\Biggl[\sqrt{\frac{1}{9}
+\frac{n}{2U_n^2[\log(3J_n(J_n+1))+t]}}-\frac{1}{3}\Biggr]\Biggr\}.\nonumber
\end{eqnarray}

Note that we allow $U_n, \eta_{1,n}, \eta_{2,n}$ and $\beta_n^{-1}$
to increase as $n$ increases. However, if those quantities are
small, the upper bound in (\ref{eqn:peoracle}) will be tighter.
\begin{theorem} \label{thm:peoracle}
Suppose Assumptions \ref{apn:errorterm} and \ref{apn:basis} hold.
For any given $0\leq\gamma<1/2$, $\eta_{2,n}>0$, $\rho\geq0$ and
$t>0$, let $\bolds{\hat\theta}_n$ be the $l_1$-PLS estimator
defined in (\ref{eqn:thetahat}) with tuning parameter
%
%e6.5 ###
\begin{eqnarray}\label{eqn:lambdacondition}
\lambda_n&\geq&\frac{8\max\{3c,2(\eta_{1,n}+\eta_{2,n})\}U_n(\log
6J_n+t)}
{(1-2\gamma)n}\nonumber\\[-8pt]\\[-8pt]
&&{}+\frac{12\max\{\sigma,(\eta_{1,n}+\eta_{2,n})\}}{(1-2\gamma)}
\sqrt{\frac{2(\log6J_n+t)}{n}}.\nonumber
\end{eqnarray}
Suppose Assumption \ref{apn:grammatrix} holds with $\rho\beta_n\leq
1$. Let $\Theta_n$ be the set defined in (\ref{eqn:thetan2}) and
assume $\Theta_n$ is nonempty.
If
%
%e6.6 ###
\begin{equation}\label{eqn:jncondition}
\frac{\log2J_n}{n}\leq
\frac{2(1-2\gamma)^2}{27U_n^2-10\gamma-22},
\end{equation}
then with probability at least $1-\exp(-k'_nn)-\exp(-t)$, we have
%
%e6.7 ###
\begin{equation}\label{eqn:peoracle}
L(\Phi_n\bolds{\hat\theta}_n)\leq
\min_{\bolds{\theta}\in\Theta_n}\biggl[L(\Phi_n\bolds{\theta})
+K_n
\frac{|M_{\rho\lambda_n}(\bolds{\theta})|}{\beta_n}\lambda_n^2\biggr],
\end{equation}
where $k_n'=13(1-2\gamma)^2/[6(27U_n^2-10\gamma-22)]$ and
$K_n=[40\gamma(12\beta_n\rho+2\gamma+5)]/[(1-2\gamma)(2\gamma+19)]
+130(12\beta_n\rho+2\gamma+5)^2/[9(2\gamma+19)^2]$.

Furthermore, suppose $E(\Phi_n^{(2)}(X,A)^T|X)=\mathbf{0}$ a.s. If
Assumption \ref{apn:grammatrix_Topt} holds with
the same $\beta_n$ as that in Assumption \ref{apn:grammatrix}, then
with probability at least
$1-\exp(-k'_nn)-\exp(-t)$, we have
\[
E\bigl(\Phi_n^{(2)}\bolds{\hat\theta}{}^{(2)}_n-T_0\bigr)^2 \leq
\min_{\bolds{\theta}\in\Theta_n}\biggl[E\bigl(\Phi_n^{(2)}\bolds{\theta}^{(2)}-T_0\bigr)^2
+ K'_n\frac{|M_{\rho\lambda_n}^{(2)}
(\bolds{\theta})|}{\beta_n}\lambda_n^2\biggr],
\]
where $K_n'=
20(12\beta_n\rho+2\gamma+5)\{\gamma/[(1-2\gamma)(7-6\beta_n\rho)]
+[3(1-2\gamma)\beta_n\rho+10(2\gamma+5)]/[9(2\gamma+19)^2]\}$.
\end{theorem}
\begin{Remarks*}
\begin{longlist}[(1)]
\item[(1)]
Note that $K_n$ is upper bounded by a constant under the assumption
$\beta_n\rho\leq1$. In the asymptotic setting when
$n\rightarrow\infty$ and $J_n\rightarrow\infty$,
(\ref{eqn:peoracle}) implies that
$L(\Phi_n\bolds{\hat\theta}_n)-\min_{\bolds{\theta}\in\mathbb
{R}^{J_n}}L(\Phi_n\bolds{\theta})\rightarrow^p
0$ if (i)
$|M_{\rho\lambda_n}(\bolds{\theta}^o)|\lambda_n^2/\beta_n=o(1)$,
(ii) $U_n^2\log J_n/n\leq k_1$ and
$|M_{\rho\lambda_n}(\bolds{\theta}^o)|\leq
k_2\beta_n\sqrt{n/(U_n^2\log J_n)}$ for some sufficiently small
positive constants $k_1$ and $k_2$ and (iii)
$\lambda_n\geq k_3\max\{1,\eta_{1,n}+\eta_{2,n}\}\sqrt{\log J_n/n}$ for
a sufficiently large
constant $k_3$, where $\bolds{\theta}^o\in[\bolds{\theta}^*_n]$ (take $t=\log J_n$).

\item[(2)] Below we briefly discuss Assumptions \ref{apn:basis}--\ref
{apn:grammatrix_Topt}.

Assumption \ref{apn:basis} is very similar to assumption \hyperlink{ap:basis}{(2)}
in Theorem \ref{thm:finaloracle} (which is used to prove the
concentration of the sample mean around the true mean), except that
$U_n$ and $\eta_{1,n}$ may increase as $n$ increases.
This relaxation allows the use of basis functions for which the sup
norm $\max_j\|\phi_j\|_\infty$ is increasing in $n$ [e.g., the wavelet
basis used in example (4) of the simulation studies].

Assumption \ref{apn:grammatrix} is a generalization of condition (\ref
{eqn:grammatrix1}) [which has been discussed in remark (4) following
Theorem \ref{thm:finaloracle}]
to the case where $J_n$ may increase in $n$ and the sparsity of a
parameter is measured by the number of ``large'' components as
described at the beginning of this section. This condition is used to
avoid the collinearity problem. It is easy to see that when $\rho=0$
and $\beta_n$ is fixed in $n$, this assumption simplifies to
condition~(\ref{eqn:grammatrix1}).

Assumption \ref{apn:grammatrix_Topt} puts a strengthened constraint on
the linear model of the treatment effect part, as compared to
Assumption \ref{apn:grammatrix}.
This assumption, together with Assumption \ref{apn:grammatrix}, is
needed in deriving the upper bound for the mean square error of the
treatment effect estimator. It is easy to verify that if $E[\Phi_n^T\Phi
_n]$ is positive definite, then both Assumptions \ref{apn:grammatrix} and \ref
{apn:grammatrix_Topt} hold.
Although the result is about the treatment effect part, which is
asymptotically independent of the main effect of $X$ (when $E[\Phi
_n^{(2)}(X,A)|X]=\mathbf{0}$ a.s.), we still need Assumption \ref
{apn:grammatrix} to show that the cross product term
$E_n[(\Phi_n^{(1)}\bolds{\hat\theta}{}^{(1)}_n-\Phi
^{(1)}_n\bolds{\theta}^{(1)})
(\Phi^{(2)}_n\bolds{\hat\theta}{}^{(2)}_n-\Phi^{(2)}_n\bolds{\theta}^{(2)})]$
is upper bounded by a quantity converging to $0$ at the desired
rate. We may use a really poor model for the main effect part
$E(Q_0(X,A)|X)$ (e.g., $\Phi^{(1)}_n\equiv1$), and Assumption~\ref{apn:grammatrix_Topt}
implies Assumption \ref{apn:grammatrix} when $\rho
=0$. This poor model only effects the constants involved in the result.
When the sample size is large (so that $\lambda_n$ is
small), the estimated ITR will be of high quality as long
as $T_0$ is well approximated.
%$\bolds{\theta}$ and
%$\bolds{\theta}^\prime\in\mathbb{R}^J$,
%E[\Phi^{(1)}(\bolds{\theta}^{(1)\prime}-\mathbf{
%||\bolds{\theta}^{(1)}||_0\geq\beta(\sum_{j\in
%M_0^{(1)}(\bolds{\theta})}\sigma_j|\theta_j^\prime-\theta_j|
%)^2 \nonumber
%E[\Phi^{(2)}(\bolds{\theta}^{(2)\prime}-\mathbf{
%||\bolds{\theta}^{(2)}||_0\geq\beta(\sum_{j\in
%M_0^{(2)}(\bolds{\theta})}\sigma_j|\theta_j^\prime-\theta_j|
%)^2.
% \label{ap:grammatrix}
\end{longlist}
\end{Remarks*}
\begin{pf*}{Proof of Theorem \ref{thm:peoracle}}
For any $\bolds{\theta}\in\Theta_n$, define the events
\begin{eqnarray*}
\Omega_1&=&\bigcap_{j=1}^{J_n}\biggl\{\frac{2(1+\gamma)}{3}\sigma_j\leq
\hat\sigma_j\leq\frac{2(2-\gamma)}{3}\sigma_j\biggr\} \qquad\mbox{[where }\hat
\sigma_j\triangleq(E_n\phi_j^2)^{1/2}\mbox{]},\\
\Omega_2(\bolds{\theta})&=&
\biggl\{\max_{j,k=1,\ldots,{J_n}}\biggl|(E-E_n)\biggl(\frac{\phi_j\phi
_k}{\sigma_j\sigma_k}\biggr)\biggr|
\leq\frac{(1-2\gamma)^2\beta_n}{120|M_{\rho\lambda_n}(\bolds{\theta})|}\biggr\},\\[-1pt]
\Omega_3(\bolds{\theta})&=&
\biggl\{\max_{j=1,\ldots,{J_n}}\biggl|E_n\biggl[(R-\Phi_n\bolds{\theta})\frac{\phi_j}{\sigma_j}\biggr]\biggr|\leq
\frac{4\gamma+1}{6}\lambda_n\biggr\}.
\end{eqnarray*}
Then there exists a
$\bolds{\theta}^o\in[\bolds{\theta}^*_n]$ such that
\begin{eqnarray*}
L(\Phi_n\bolds{\hat\theta}_n)
&=&L(\Phi_n\bolds{\theta})
+2E[(\Phi_n\bolds{\theta}^o-\Phi_n\bolds{\theta})\Phi
_n(\bolds{\theta}-\bolds{\hat\theta}_n)]
+E[\Phi_n(\bolds{\hat\theta}_n-\bolds{\theta})]^2\\[-1pt]
&\leq& L(\Phi_n\bolds{\theta})+2\max_{j=1,\ldots,J_n}\biggl|
E\biggl[\Phi_n(\bolds{\theta}^o-\bolds{\theta})\frac{\phi
_j}{\sigma_j}\biggr]\biggr|
\Biggl(\sum_{j=1}^{J_n}\sigma_j|\hat\theta_{n,j}-\theta_j|\Biggr)\\[-1pt]
&&{}+
E[\Phi_n(\bolds{\hat\theta}_n-\bolds{\theta})]^2\\[-2pt]
&\leq& L(\Phi_n\bolds{\theta})+2\gamma\lambda_n
\Biggl(\sum_{j=1}^{J_n}\sigma_j|\hat\theta_{n,j}-\theta_j|\Biggr)+
E[\Phi_n(\bolds{\hat\theta}_n-\bolds{\theta})]^2,
\end{eqnarray*}
where the first equality follows from the fact that
$E[(R-\Phi_n\bolds{\theta}^o)\phi_j]=0$ for any
$\bolds{\theta}^o\in[\bolds{\theta}^*_n]$ for $j=1,\ldots,
J_n$ and the last inequality follows from the definition
of~$\Theta_n^o$.

Based on Lemma \ref{lemma:oracle} below, we have that on the event
$\Omega_1\cap\Omega_2(\bolds{\theta})\cap\Omega_3(\bolds{\theta})$,
\[
L(\Phi_n\bolds{\hat\theta}_n) \leq
L(\Phi_n\bolds{\theta})+K_n\frac{|M_{\rho\lambda_n}
(\bolds{\theta})|}{\beta_n}\lambda_n^2.\vspace{-2pt}
\]

Similarly, when $E[\Phi^{(2)}_2(X,A)^T|X]=\mathbf{0}$, by Lemma \ref
{lemma:oracle1}, we have that on the event
$\Omega_1\cap\Omega_2(\bolds{\theta})\cap\Omega_3(\bolds{\theta})$,\vspace{-3pt}
\begin{eqnarray*}
E\bigl(\Phi_n^{(2)}\bolds{\hat\theta}{}^{(2)}_n-T_0\bigr)^2
&\leq& E\bigl(\Phi_n^{(2)}\bolds{\theta}^{(2)}-T_0\bigr)^2
+2\gamma\lambda_n
\Biggl(\sum_{j=J_n^{(1)}+1}^{J_n}\sigma_j|\hat\theta_{n,j}-\theta_j|\Biggr)\\[-1pt]
&&{} +
E\bigl[\Phi_n^{(2)}\bigl(\bolds{\hat\theta}{}^{(2)}_n-\bolds{\theta}^{(2)}\bigr)\bigr]^2\\[-1pt]
&\leq&
E\bigl(\Phi_n^{(2)}\bolds{\theta}^{(2)}-T_0\bigr)^2+K'_n\frac{|M_{\rho\lambda
_n}^{(2)}
(\bolds{\theta})|}{\beta_n}\lambda_n^2.
\end{eqnarray*}

The conclusion of the theorem follows from the union probability
bounds of the events $\Omega_1$, $\Omega_2(\bolds{\theta})$ and
$\Omega_3(\bolds{\theta})$ provided in Lemmas
\ref{lemma:omega1}, \ref{lemma:omega2} and \ref{lemma:omega3}.\vspace{-5pt}
\end{pf*}

Below we state the lemmas used in the proof of Theorem \ref{thm:peoracle}.
The proofs of the lemmas are given in Section S.4 of the supplemental
article \cite{supplement}.
\begin{lemma}\label{lemma:oracle}
Suppose Assumption \ref{apn:grammatrix} holds with $\rho\beta_n\leq
1$. Then for any $\bolds{\theta}\in\Theta_n$, on the event
$\Omega_1\cap\Omega_2(\bolds{\theta})\cap\Omega_3(\bolds{\theta})$,
we have
%
%e6.8 ###
\begin{equation}\label{eqn:hatthetabound}
\sum_{j=1}^{J_n}\sigma_j|\hat\theta_{n,j}-\theta_j|\leq
\frac{20(12\rho\beta_n+2\gamma+5)}{(1-2\gamma)(19+2\gamma)\beta_n}
|M_{\rho\lambda_n}(\bolds{\theta})|\lambda_n
\end{equation}
and
%
%e6.9 ###
\begin{equation}\label{eqn:thetariskbound}
E[\Phi_n(\bolds{\hat\theta}_n-\bolds{\theta})]^2\leq
\frac{130(12\rho\beta_n+2\gamma+5)^2
}{9(19+2\gamma)^2\beta_n}|M_{\rho\lambda_n}(\bolds{\theta})|\lambda_n^2
\end{equation}
\end{lemma}
\begin{Remark*}
This lemma implies that $\bolds{\hat\theta}_n$ is close to each
$\bolds{\theta}\in\Theta_n$ on the event
$\Omega_1\cap\Omega_2(\bolds{\theta})\cap\Omega_3(\bolds{\theta})$.
The intuition is as follows. Since $\bolds{\hat\theta}_n$
minimizes~(\ref{eqn:thetahat}), the first order conditions imply
that
$\max_j|E_n(R-\Phi_n\bolds{\hat\theta}_n)\phi_j/\hat\sigma_j|\leq
\lambda_n/2$.
Similar property holds for $\bolds{\theta}$ on the event
$\Omega_1\cap\Omega_3(\bolds{\theta})$. Assumption
\ref{apn:grammatrix} together with event
$\Omega_2(\bolds{\theta})$ ensures that there is no
collinearity in the $n\times J_n$ design matrix
$(\Phi_n(X_i,A_i))_{i=1}^n$. These two aspects guarantee the
closeness of $\bolds{\hat\theta}_n$ to~$\bolds{\theta}$.
\end{Remark*}
\begin{lemma}\label{lemma:oracle1}
Suppose $E[\Phi_n^{(2)}(X,A)^T|X]=\mathbf{0}$ a.s. and
Assumptions \ref{apn:grammatrix} and \ref{apn:grammatrix_Topt} hold
with $\rho\beta_n\leq1$.
Then for any $\bolds{\theta}\in\Theta_n$, on the event
$\Omega_1\cap\Omega_2(\bolds{\theta})\cap\Omega_3(\bolds{\theta})$,
we have
%
%e6.10 ###
\begin{equation}
\label{eqn:hatthetabound1}
\sum_{j=J_n^{(1)}+1}^{J_n}\sigma_j|\hat\theta_{n,j}-\theta_j|\leq
\frac{10(12\beta_n\rho+2\gamma+5)}{(1-2\gamma)(7-6\beta_n\rho)\beta
_n}\bigl|M_{\rho\lambda_n}^{(2)}
(\bolds{\theta})\bigr|\lambda_n
\end{equation}
and
%
%e6.11 ###
\begin{eqnarray}\label{eqn:thetariskbound1}
\hspace*{28pt}&&E\bigl[\Phi_n^{(2)}\bigl(\bolds{\hat\theta}{}^{(2)}_n-\bolds{\theta}^{(2)}\bigr)\bigr]^2
\nonumber\\[-8pt]\\[-8pt]
\hspace*{28pt}&&\qquad\leq
\frac{20(12\rho\beta_n+2\gamma+5)[3(1-2\gamma)\beta_n\rho+10(2\gamma+5)]
}{9(2\gamma+19)^2\beta_n}\bigl|M_{\rho\lambda_n}^{(2)}(\bolds{\theta})\bigr|\lambda_n^2.\nonumber
\end{eqnarray}
\end{lemma}
\begin{lemma}\label{lemma:omega1}
Suppose Assumption \ref{apn:basis}\textup{(a)} and inequality
(\ref{eqn:jncondition}) hold. Then $\mathbf{P}(\Omega_1^C)\leq\exp
(-k'_nn)$, where $k'_n = 13(1-2\gamma)^2/[6(27U_n^2-10\gamma-22)]$.
\end{lemma}
\begin{lemma}\label{lemma:omega2}
Suppose Assumption \ref{apn:basis}\textup{(a)} holds. Then for any
$t>0$ and $\bolds{\theta}\in\Theta_n$,
$\mathbf{P}(\{\Omega_2(\bolds{\theta})\}^C)\leq\exp(-t)/3$.
\end{lemma}
\begin{lemma}\label{lemma:omega3}
Suppose Assumptions \ref{apn:errorterm} and \ref{apn:basis} hold.
For any $t>0$, if $\lambda_n$
satisfies condition (\ref{eqn:lambdacondition}), then for any
$\bolds{\theta}\in\Theta_n$, we have
$\mathbf{P}(\{\Omega_3(\bolds{\theta})\}^C)\leq2\exp(-t)/3$.
\end{lemma}

\subsection{\texorpdfstring{Design of simulations in Section \lowercase{\protect\ref{sec:simulation}}}%
{Design of simulations in Section 5.1}}
\label{sec:simdesign}

In this section, we present the detailed simulation design of the
examples used in Section \ref{sec:simulation}. These examples
satisfy all assumptions listed in the theorems [it is easy to verify
that for examples (1)--(3). Validity of the assumptions for
example (4) is addressed in the remark after example (4)]. In addition,
$\Theta_n$ defined in (\ref{eqn:oracleset}) is nonempty as long as
$n$ is sufficiently large (note that the constants involved in
$\Theta_n$ can be improved and are not that meaningful. We focused
on a presentable result instead of finding the best constants).

In examples (1)--(3),
$X=(X_1,\ldots,X_5)$ is uniformly distributed on $[-1,1]^5$. The
treatment $A$ is then generated independently of $X$ uniformly from $\{
-1,1\}$. Given $X$ and $A$, the response $R$ is
generated from a normal distribution with mean
$Q_0(X,A)=1+2X_1+X_2+0.5X_3 + T_0(X,A)$ and variance~$1$. We consider
the following three examples for
$T_0$:
\begin{longlist}[(1)]
\item[(1)]$T_0(X,A) = 0$ (i.e., there is no treatment
effect).
\item[(2)] $T_0(X,A) = 0.424(1-X_1-X_2)A$.
\item[(3)]$T_0(X,A) = 0.446 \operatorname{sign}(X_1)(1-X_1)^2A$.
\end{longlist}
Note that in each example $T_0(X,A)$ is equal to the treatment effect
term, $ Q_0(X,A)-E[Q_0(X,A)|X]$. We approximate $Q_0$ by $\mathcal{Q}=\{
(1, X, A, XA)\bolds{\theta}\dvtx\break\bolds{\theta}\in\mathbb{R}^{12}\}$.
Thus, in examples (1) and (2) the treatment effect term $T_0$ is correctly modeled,
while in example (3) the treatment effect term $T_0$ is misspecified.

The parameters in examples (2) and (3) are chosen to reflect a medium
effect size according to Cohen's d index. When there are two
treatments, the Cohen's d effect size index is defined as the
standardized
difference in mean responses between two treatment groups, that is,
\[
\mathrm{es}=\frac{E(R|A=1)-E(R|A=-1)}{([\operatorname{Var}(R|A=1)+\operatorname{Var}(R|A=-1)]/2)^{1/2}}.
\]
Cohen \cite{cohen1988} tentatively defined the effect size as ``small''
if the Cohen's d index is $0.2$, ``medium'' if the index is $0.5$
and ``large'' if the index is $0.8$.

In example (4), $X$ is uniformly distributed on $[0,1]$. Treatment $A$
is generated independently of $X$ uniformly from $\{-1,1\}$. The
response $R$ is generated from a normal distribution
with mean $Q_0(X,A)$ and variance $1$, where $Q_0(X,1)
=\sum_{j=1}^8\vartheta_{(1),j}1_{X<u_{(1),j}}$, $Q_0(X,-1)
=\sum_{j=1}^8\vartheta_{(-1),j}1_{X<u_{(-1),j}}$, and $\vartheta$'s
and $u$'s are parameters specified in (\ref{eqn:parameter}). The
effect size is small:
%&\vartheta_{(0),1}=-0.4260, \vartheta_{(0),2}= 2.8856,
%&\vartheta_{(1),1}=2.0822, \vartheta_{(1),2}=-0.7318,
%&u_{(0),1}=0.1408, u_{(0),2}=0.9902, u_{(0),3}=0.2807,
%u_{(0),4}=0.4929, u_{(0),5}=0.4651;\nonumber\\
%& u_{(1),1}=0.9934, u_{(1),2}=0.1191, u_{(1),3}=0.2509,
%u_{(1),4}=0.7541, u_{(1),5}=0.6660.
%
%e6.12 ###
\begin{eqnarray}\label{eqn:parameter}\qquad
&&\bigl(\vartheta_{(1),1},\ldots, \vartheta_{(1),8}\bigr)\nonumber\\
&&\qquad= (-0.781, 0.730,
0.635, 0.512, -2.278, 1.347, 1.155,
-0.030); \nonumber\\
&&\bigl(\vartheta_{(-1),1},\ldots, \vartheta_{(-1),8}\bigr)\nonumber\\
&&\qquad= (-2.068, 1.520,
-0.072,\nonumber\\
&&\hspace*{37.7pt} -0.637, 1.003, -0.611, -0.305, 1.016); \\
&&\bigl(u_{(1),1},\ldots, u_{(1),8}\bigr)\nonumber\\
&&\qquad= (0.028, 0.144, 0.171, 0.298,
0.421, 0.443, 0.463, 0.758); \nonumber\\
&&\bigl(u_{(-1),1},\ldots, u_{(-1),8}\bigr)\nonumber\\
&&\qquad= (0.061, 0.215, 0.492, 0.544,
0.6302, 0.650, 0.785, 0.909).\nonumber
\end{eqnarray}

We approximate $Q_0$ by Haar wavelets,
\[
\theta_{(0),0}h_0(X)+\sum_{lk}\theta_{(0),lk}h_{lk}(X)+
\biggl(\theta_{(0),1}h_0(X)+\sum_{lk}\theta_{(1),lk}h_{lk}(X)\biggr)A,
\]
where $h_0(x)=1_{x\in[0,1]}$ and
$h_{lk}(x)=2^{l/2}(1_{2^lx\in[k+1/2,k+1)}-
1_{2^lx\in[k,k+1/2)})$ for $l=0,\ldots,\bar l_n$, and $\theta
_{(\cdot),\cdot}\in\mathbb{R}$ are parameters. We choose
$\bar l_n =\lfloor3\log_2 n/4\rfloor-2$. For a given $l$ and sample
$(X_i,A_i,R_i)_{i=1}^n$, $k$ takes integer values from $\lfloor
2^l\min_iX_i\rfloor$ to $\lceil2^l\max_iX_i\rceil-1$. Then $J_n =
2^{\lfloor3\log_2n/4\rfloor}\leq n^{3/4}$.

\begin{Remark*}
In example (4), we allow the number of basis functions $J_n$ to
increase with $n$. The corresponding theoretical result can be
obtained by combining Theorems \ref{thm:bound2} and
\ref{thm:peoracle}. Below we demonstrate the validation of the
assumptions used in the theorems.

Theorem \ref{thm:bound2} requires that the randomization probability
$p(a|x)\geq S^{-1}$ for a positive constant for all $(x,a)$ pairs
and the margin condition (\ref{eqn:noise}) or (\ref{eqn:noise2})
holds. According the generative model, we have that $p(a|x)=1/2$ and
condition (\ref{eqn:noise2}) holds.

Theorem \ref{thm:peoracle} requires Assumptions \ref
{apn:errorterm}--\ref{apn:grammatrix_Topt} hold and $\Theta
_n$ defined
in (\ref{eqn:thetan2}) is nonempty. Since we consider normal error
terms, Assumption \ref{apn:errorterm} holds. Note that the basis
functions used in Haar wavelet are orthogonal. It is also easy to
verify that Assumptions \ref{apn:grammatrix} and
\ref{apn:grammatrix_Topt} hold with $\beta_n=1$ and Assumption
\ref{apn:basis} holds with $U_n = n^{3/8}/2$ and $\eta_{1,n}\leq
\mbox{constant} + \mbox{constant}\times\|\bolds{\theta}^*_n\|_0$ [since each
$|\phi_j\bolds{\theta}^*_{n,j}| =|\phi_j E(\phi_jR)|\leq
\mbox{constant}\times|\phi_j| E|\phi_j|\leq O(1)$]. Since $Q_0$ is piece-wise
constant, we can also verify
that $\|\bolds{\theta}^*_n\|_0\leq O(\log n)$. Thus, for sufficiently
large $n$, $\Theta_n$
is nonempty and (\ref{eqn:jncondition}) holds. The RHS of
(\ref{eqn:lambdacondition}) converges to zero as
$n\rightarrow\infty$.
\end{Remark*}
\end{appendix}

\section*{Acknowledgments}
The authors thank Martin Keller and the investigators of the
Nefazodone-CBASP trial for use of their data. The authors also thank
John Rush, MD, for the technical support and Bristol-Myers Squibb
for helping fund the trial. The authors acknowledge the help of
the reviewers and of Eric B. Laber and Peng Zhang in improving this paper.

% AOS,AOAS: If there are supplements please fill:
%
\begin{supplement}
\stitle{Supplement to ``Performance guarantees for individualized
treatment rules''}
\slink[doi]{10.1214/10-AOS864SUPP}
\sdatatype{.pdf}
\sfilename{aos864\_suppl.pdf}
\sdescription{This supplement contains four sections. Section S.1
discusses the problem with over-fitting due to the potentially
large number of pretreatment variables (and/or complex approximation
space for $Q_0$) mentioned in Section \ref{sec:lasso}. Section S.2
provides modifications of
the $l_1$-PLS estimator $\bolds{\hat\theta}_n$ when some
coefficients are not penalized and discusses how to obtain results
similar to inequality (\ref{eqn:peoracle}) in this case. Section S.3
provides extra four simulation examples
based on data from the Nefazodone-CBASP
trial \cite{keller2000}. Section~S.4 provides proofs of
Lemmas~\ref{lemma:oracle}--\ref{lemma:omega3}.}
\end{supplement}

%suskaldyti doi

% imsref loaded by lrinkeviciute, 2011-02-17 08:45:23
% imsref loaded by lrinkeviciute, 2011-02-17 08:49:12
%

%
\printaddresses

\end{document}